\RequirePackage{etoolbox}
\csdef{input@path}{{style/}{graphics/}}
\documentclass[aos,MSNbibl,nameyear,secthm,dvips]{arximspdf}
\usepackage{mathbh}
\usepackage{graphicx}

\doi{10.1214/14-AOS1304}% Updated by VTEXPTS2LaTeX.exe, 23.02.2015 08:41
\volume{43}
\issue{2}
\pubyear{2015}
\firstpage{878}
\lastpage{902}
\docsubty{FLA}

\makeatletter
\newproclaim{rem}{Remark}
\newcommand{\ds}{\displaystyle}

\newtheorem{prop}[thm]{Proposition}
\newtheorem{cor}[thm]{Corollary}
\newcommand{\cov}{\operatorname{Cov}}

\newcommand{\dt}[1]{\dot{#1}}
\newcommand{\dd}{\mathrm{d}}
\newcommand{\iint}{\int\!\!\!\int}
\makeatother

\begin{document}
\begin{frontmatter}

\title{Asymptotically distribution-free goodness-of-fit testing for
tail copulas}
\runtitle{Distribution-free GOF testing for tail copulas}

\begin{aug}
% Corresponding author: John Einmahl - j.h.j.einmahl@tilburguniversity.edu% Updated by VTEXPTS2LaTeX.exe, 24.02.2015 13:27
%j.h.j.einmahl@tilburguniversity.edu% Updated by VTEXPTS2LaTeX.exe,
%23.02.2015 08:41
\author[A]{\fnms{Sami Umut}~\snm{Can}\thanksref{m11,M1}\ead[label=e1]{s.u.can@uva.nl}},
\author[B]{\fnms{John H. J.}~\snm{Einmahl}\thanksref{m2}\corref{}\ead[label=e2]{j.h.j.einmahl@tilburguniversity.edu}},
\author[C]{\fnms{Estate~V.}~\snm{Khmaladze}\thanksref{m3}\ead[label=e3]{estate.khmaladze@vuw.ac.nz}}
\and
\author[A]{\fnms{Roger J. A.}~\snm{Laeven}\thanksref{m11}\ead[label=e4]{r.j.a.laeven@uva.nl}}
\runauthor{Can, Einmahl, Khmaladze and Laeven}
\affiliation{University of Amsterdam\thanksmark{m11}, Tilburg University\thanksmark{m2}\\ and  Victoria
University of Wellington\thanksmark{m3}}
\address[A]{S. U. Can\\
R. J. A. Laeven\\
Faculty of Economics \& Business\\
Section Actuarial Science\\
University of Amsterdam\\
Valckenierstraat 65\\
1018 XE Amsterdam\\
The Netherlands\\
\printead{e1}\\
\phantom{E-mail:\ }\printead*{e4}}
\address[B]{J. H. J. Einmahl\\
Department of Econometrics \& OR\\
\quad and CentER\\
Tilburg University\\
P.O. Box 90153\\
5000 LE Tilburg\\
The Netherlands\\
\printead{e2}}
\address[C]{E. V. Khmaladze\\
School of Mathematics, Statistics \& OR\\
Victoria University of Wellington\\
P.O. Box 600\\
Wellington\\
New Zealand\\
\printead{e3}}
\end{aug}
\thankstext{M1}{Supported in part by the Netherlands Organization for
Scientific Research under grant NWO Vidi 2009.}

% HISTORY:
%
\received{\smonth{11} \syear{2013}}% Updated by VTEXPTS2LaTeX.exe,
%23.02.2015 08:41
%
\revised{\smonth{12} \syear{2014}}% Updated by VTEXPTS2LaTeX.exe,
%23.02.2015 08:41

% ABSTRACT
\begin{abstract}
Let $(X_1,Y_1), \ldots, (X_n,Y_n)$ be an i.i.d. sample from a
bivariate distribution function that lies in the max-domain of
attraction of an extreme value distribution. The asymptotic joint
distribution of the standardized component-wise maxima $\bigvee_{i=1}^n
X_i$ and $\bigvee_{i=1}^n Y_i$ is then characterized by the marginal
extreme value indices and the tail copula $R$. We propose a procedure
for constructing asymptotically distribution-free goodness-of-fit tests
for the tail copula $R$. The procedure is based on a transformation of
a suitable empirical process derived from a semi-parametric estimator
of $R$. The transformed empirical process converges weakly to a
standard Wiener process, paving the way for a multitude of
asymptotically distribution-free goodness-of-fit tests. We also extend
our results to the $m$-variate ($m>2$) case. In a simulation study we
show that the limit theorems provide good approximations for finite
samples and that tests based on the transformed empirical process have
high power.
\end{abstract}

% KEYWORDS
\end{frontmatter}

%s1 #&#
\section{Introduction}\label{sec.intro}

Let $(X_1,Y_1), \ldots, (X_n,Y_n)$ be an i.i.d. sample from a
bivariate distribution function (d.f.) $F$ with marginal d.f.'s $F_1(x) =
F(x,\infty)$ and $F_2(y)=F(\infty,y)$ for $x, y \in\mathbb{R}$.
Suppose that $F$ is in the max-domain of attraction of some bivariate
d.f. $G$ with nondegenerate marginals. That is, suppose that there exist
normalizing sequences $a_1(n), a_2(n)>0$ and $b_1(n), b_2(n) \in\mathbb
{R}$ such that
%
%e1 #&#
\begin{equation}
\label{biv.doa}
P \biggl( \frac{\bigvee_{i=1}^n X_i - b_1(n)}{a_1(n)} \le x, \frac
{\bigvee_{i=1}^n Y_i - b_2(n)}{a_2(n)} \le y \biggr)
\rightarrow G(x,y)
\end{equation}
as $n \to\infty$, for all continuity points $(x,y) \in\mathbb{R}^2$
of $G$. Of course, (\ref{biv.doa}) is equivalent to
%
%e2 #&#
\begin{equation}
\label{biv.doa.2} F^n \bigl(a_1(n)x + b_1(n),
a_2(n)y + b_2(n) \bigr) \rightarrow G(x,y),
\end{equation}
and the d.f. $G$ is, by definition, an extreme value d.f.

It is a classical result in extreme value theory [see \citet{dehaanferreira2006}, Theorem~1.1.3] that the normalizing sequences
$a_1, b_1$ and $a_2, b_2$ can be chosen in such a way that the marginal
d.f.'s $G_1(x) = G(x,\infty)$ and $G_2(y) = G(\infty,y)$ are of the form
%
%e3 #&#
\begin{eqnarray}
G_1(x) &= &\exp\bigl\{-(1+\gamma_1x)^{-1/\gamma_1}
\bigr\}, \qquad 1+\gamma_1x > 0,
\nonumber
\\[-8pt]
\label{hazaza}
\\[-8pt]
\nonumber
G_2(y) &=& \exp\bigl\{-(1+\gamma_2y)^{-1/\gamma_2}
\bigr\}, \qquad 1+\gamma_2y > 0
\end{eqnarray}
for\vspace*{1pt} some $\gamma_1, \gamma_2 \in\mathbb{R}$. [Here, and in the rest of
the paper, expressions of the form $(1 + \gamma\,\cdot\,)^{1/\gamma}$
should be interpreted as $\exp(\cdot)$ when $\gamma= 0$.] We will
assume throughout that the normalizing sequences are chosen in this
way. Then $G$ is necessarily continuous, as it has continuous marginal
d.f.'s, and the equivalent convergences (\ref{biv.doa}) and (\ref{biv.doa.2}) hold for all $(x,y) \in[-\infty,\infty]^2$. Also, $G$ can
be fully characterized by the marginal extreme value indices $\gamma
_1$, $\gamma_2$ and a description of the dependence structure between
the marginal d.f.'s $G_1$ and $G_2$. Due to \citet{dehaanresnick1977}, it
is known that the class of possible dependence structures for bivariate
extreme value distributions does not form a finite-dimensional
parametric family. Nevertheless, there are various equivalent ways of
describing extreme value (or tail) dependence structures, each with its
own advantages in applications. For an overview, we refer to \citet{beirlantetal2004}, Chapter~8 or \citet{dehaanferreira2006}, Part~II.

In this paper, we will focus on one possible description of the
bivariate tail dependence structure, namely the \textit{tail copula}. For
a bivariate extreme value d.f. $G$ with marginal d.f.'s as given in (\ref
{hazaza}), the tail copula $R$ is defined as
%
%e4 #&#
\begin{equation}
\label{R.def}
R(x,y) = x + y + \log G \biggl( \frac{x^{-\gamma_1} - 1}{\gamma_1},
\frac
{y^{-\gamma_2} - 1}{\gamma_2} \biggr), \qquad (x,y) \in[0, \infty)^2.
\end{equation}
We say that a bivariate d.f. $F$ belonging to the domain of attraction of
$G$ has associated tail copula $R$. It is clear that tail copulas are
not copula functions in the usual sense (since they are not
distribution functions of probability measures, e.g.), yet they fully
capture the asymptotic dependence structure of the component-wise
maxima, just like copulas capture the dependence structure of random
vectors. Indeed, it is easily checked that $G(x,y) = C_G(G_1(x),
G_2(y))$, with
%
%e5 #&#
\begin{equation}
\label{G.copula}
C_G(u,v) = uv\exp\bigl\{R(-\log u, -\log v)\bigr\},
\qquad (u,v) \in(0,1]^2.
\end{equation}
In other words, $G$ is the unique d.f. characterized by the marginal d.f.'s
(\ref{hazaza}) and the copula (\ref{G.copula}).

We conclude that the asymptotic joint behavior of the standardized\break component-wise maxima $\bigvee_{i=1}^n X_i$ and $\bigvee_{i=1}^n Y_i$ is
fully characterized by the margi\-nal extreme value indices $\gamma_1,
\gamma_2$ appearing in (\ref{hazaza}) and the tail copula $R$ defined
in (\ref{R.def}). Statistical inference about extreme value indices is
a classical and well-studied problem in univariate extreme value
theory; we refer to \citet{beirlantetal2004}, Chapters~4 and~5 or
\citet{dehaanferreira2006}, Chapter~3 for more information. There is
also a growing literature on inference about the tail dependence
structure; see \citet{beirlantetal2004}, Chapter~9  or
\citet{dehaanferreira2006}, Chapter~7, for an overview. In this paper, we
will focus on inference about $R$. In particular, we will propose a
semi-parametric estimator of $R$, describe a transformation of the
empirical process derived from it and demonstrate how this transformed
empirical process can serve as a basis to construct asymptotically
distribution-free goodness-of-fit tests for $R$.

%s1.1 #&#
\subsection{More on tail dependence}
The tail copula $R$ can also be obtained (and its domain extended) in
the following way from the d.f. $F$:
%
%e6 #&#
%e7 #&#
\begin{eqnarray}
R(x,y) &=& \lim_{t\to\infty} tP \bigl(1-F_1(X)
<x/t, 1-F_2(Y)<y/t \bigr),
\nonumber
\\[-8pt]
\label{R.def.2}
\\[-8pt]
\eqntext{(x,y) \in[0,\infty]^2 \setminus\bigl\{(\infty,\infty)\bigr
\},} %
\end{eqnarray}
where $(X,Y)$ denotes a random vector with d.f. $F$. If $F$ has
continuous marginals, (\ref{R.def.2}) can also be written as
%
%e8 #&#
\begin{equation}
\label{R.def.3}
R(x,y) = \lim_{t \to\infty} t \overline{C}_F(x/t,
y/t), \qquad (x,y) \in[0,\infty]^2 \setminus\bigl\{(\infty,\infty)\bigr
\},
\end{equation}
where $\overline{C}_F$ denotes the ``survival copula'' of $F$, that is,
the copula associated with $(-X,-Y)$. Observe that $R(x,\infty) =
R(\infty,x) = x$  for all  $x \in[0,\infty)$ and $0 \le R(x,y)
\le x \wedge y$  for all  $(x,y) \in[0,\infty]^2\setminus\{(\infty
,\infty)\}$. It is also clear from (\ref{R.def.2}) that $R$ is
homogeneous of order 1, so the restriction of $R$ on, for example,
$[0,1]^2$ determines $R$ on its entire domain. The characterization
(\ref{R.def.2}) stems from \citet{huang1992}, where it is used to
derive a nonparametric estimator for $R$. We will use an alternative,
semi-parametric estimator better suited for our purposes; see Section~\ref{sec.estimator}.

The value $R(1,1)$ is known in the applied extreme value literature as
the (\textit{upper}) \textit{tail dependence coefficient} and is widely used as a
measure of tail dependence. When $R(1,1)=0$, which is equivalent to $R
\equiv0$ on $[0,\infty)^2$, we call $X$ and $Y$ tail independent. When
$R(1,1) > 0$, we say that $X$ and $Y$ exhibit tail dependence. Other
ways of describing the tail dependence structure include the
\textit{stable tail dependence function}, the \textit{exponent measure}, the
\emph{spectral measure} and the \emph{Pickands dependence function};
see the monographs \citet{kotznadarajah2000}, \citet{beirlantetal2004}, \citet{dehaanferreira2006} and the many
references therein.

We also note here that the function $R$ generates a $\sigma$-finite
measure, which we will also, without confusion, denote by $R$, on Borel
subsets of $[0,\infty]^2 \setminus\{(\infty,\infty)\}$, through the identity
%
%e9 #&#
\begin{equation}
\label{Lambda.def}
R\bigl([0,x] \times[0,y]\bigr) := R(x,y), \qquad (x,y) \in[0,
\infty]^2 \setminus \bigl\{(\infty,\infty)\bigr\}.
\end{equation}

%s1.2 #&#
\subsection{Goodness-of-fit testing}
In the literature and in practice, often a parametric model is used for
the tail copula $R$; see, for example, \citet{colestawn1991} or \citet{joeetal1992}. Testing the goodness-of-fit of the parametric model to
a given data sample is therefore an important problem with abundant
applications in many fields such as insurance and risk management,
finance and econometrics and hydrology and meteorology. In this paper,
we develop a procedure for constructing asymptotically
distribution-free goodness-of-fit tests for the tail copula $R$ of a
bivariate d.f. $F$. We consider null hypotheses of the form $R \in
\mathcal{R} = \{R_{\bolds{\theta}} \dvtx  \bolds{\theta} \in\Theta
\}$, where $\mathcal{R}$ is a parametric family of tail copulas. Of
course, by taking the parameter space $\Theta$ to consist of a single
point, our results can also be used to test the goodness-of-fit of a
fully specified tail copula\vspace*{1pt} to the data.

Our approach is based on a semi-parametric estimator $\widehat{R}_n$ of
$R$, to be defined below. We consider a suitably normalized difference,
$\widehat{\eta}_n$, between $\widehat{R}_n$ and $R_{\widehat{\bolds
{\theta}}}$ (with $\widehat{\bolds{\theta}}$ denoting a suitable
estimator of $\bolds{\theta}$), and we show that, under the null
hypothesis, a proper transformation of $\widehat{\eta}_n$ converges
weakly to a standard Wiener process $W$. This fundamental result allows
one to construct a myriad of goodness-of-fit tests based on comparisons
of appropriate functionals of $\widehat{\eta}_n$ (the test statistics
the practitioner may prefer to use) with the same functionals of $W$.
We emphasize that, since $W$ is a \emph{standard} Wiener process, our
approach leads to \emph{asymptotically distribution-free}
goodness-of-fit tests: under the null hypothesis, the asymptotic
distributions of the test statistics do not depend on $\mathcal{R}$ or
the true $\bolds{\theta}$. A simulation study confirms the
applicability of our approach for finite samples.

Testing (and estimation) problems for the tail copula have been studied
in the recent literature. In \citet{einmahletal2006} the existence of
$R$ is tested, rather than its membership of a parametric family. In
\citet{dehaanetal2008} a specific Cram\'er--von Mises type statistic
for $R \in \{R_{\bolds{\theta}} \dvtx  \bolds{\theta} \in\Theta\}
$ is studied for two-dimensional data and a one-dimensional parameter;
the test statistic has a complicated limiting distribution under the
null hypothesis. In \citet{einmahletal2012} it is assumed that $R \in
\{R_{\bolds{\theta}} \dvtx  \bolds{\theta} \in\Theta\}$, and it
is then tested if $R$ is a member of a smaller parametric family,
obtained by setting some components of $\bolds{\theta}$ equal to
fixed values.

The remainder of the paper is organized as follows. In Section~\ref{sec.estimator}, we
describe the semi-parametric estimator $\widehat{R}_n$, introduce the
empirical process $\widehat{\eta}_n$, which is the normalized
difference between $\widehat{R}_n$ and $R_{\widehat{\bolds{\theta
}}}$, and describe the weak limit $\widehat{\eta}$ of $\widehat{\eta
}_n$ as $n \to\infty$. In Section~\ref{sec.trans}, we describe our key
transformation from $\widehat{\eta}$ into a standard Wiener process. In
Section~\ref{sec.gof}, we show that the same transformation (or rather an empirical
version of it, with unknown parameters replaced by estimators) applied
to $\widehat{\eta}_n$ produces a process whose weak limit is a standard
Wiener process. This is our main result. In Section~\ref{sec.multi}, we extend this
result to the $m$-dimensional setting, for $m > 2$. Finally, in Section~\ref{sec.sim}, we demonstrate through Monte Carlo simulations the applicability of
our limit theorems in finite samples and the high power properties of
tests based on our results. Proofs are deferred to Section~\ref{sec.proofs}. The paper
is supplemented by an online appendix, see \citet{supp}, which
contains some details suppressed in Section~\ref{sec.estimator} as well
as technical specifics about the Monte Carlo simulations, including the
computer code. %\footnote{The appendix is accessible at

%%%%%%%%%%%%%%%%%%%%%%%%%%%%%%%%%%%%%%%%%%%%%%%%%%%%%%%%%%%%%%%%%%%%%%%%%%%%%%%%%%%%%%%%%%%%%%%%
%s2
%s2 #&#
\section{An estimator for $R$ and its asymptotic behavior}\label{sec.estimator}
%%%%%%%%%%%%%%%%%%%%%%%%%%%%%%%%%%%%%%%%%%%%%%%%%%%%%%%%%%%%%%%%%%%%%%%%%%%%%%%%%%%%%%%%%%%%%%%%

As in Section~\ref{sec.intro}, we let $(X_1,Y_1), \ldots, (X_n,Y_n)$
denote an i.i.d. sample from a bivariate d.f. $F$ with marginal d.f.'s $F_1$
and $F_2$. We assume that the bivariate domain of attraction condition
(\ref{biv.doa}) holds, with the normalizing sequences $a_1, b_1$ and
$a_2, b_2$ chosen such that the marginal d.f.'s $G_1$ and $G_2$ are as in
(\ref{hazaza}). Taking logarithms in (\ref{biv.doa.2}), and replacing
the discrete index $n$ by a continuous index $t>0$, we obtain
\[
\lim_{t \to\infty} t\bigl[1 - F\bigl(a_1(t)x +
b_1(t), a_2(t)y + b_2(t)\bigr)\bigr] = -\log
G(x,y), \qquad (x,y) \in\mathbb{R}^2.
\]
Combining this with the corresponding marginal results and (\ref
{G.copula}) leads to
\begin{eqnarray*}
&& \lim_{t \to\infty} t P \bigl( X_1 > a_1(t)x
+ b_1(t), Y_1 > a_2(t)y + b_2(t)
\bigr)
\\
&&\qquad= R\bigl(-\log G_1(x), -\log G_2(y)\bigr),
\end{eqnarray*}
or equivalently,
\[
\lim_{t \to\infty} t P \bigl( \widetilde{X}_1(t) \leq x,
\widetilde {Y}_1(t) \leq y \bigr) = R(x,y),
\]
with
%
%e10 #&#
\begin{eqnarray}
\widetilde{X}_i(t) &=& \biggl[ \biggl( 1 +
\gamma_1 \frac{X_i -
b_1(t)}{a_1(t)} \biggr)\vee0 \biggr]^{-1/\gamma_1},
\nonumber
\\[-8pt]
\label{tildes}
\\[-8pt]
\nonumber
\widetilde{Y}_i(t) &=& \biggl[ \biggl( 1 +
\gamma_2 \frac{Y_i -
b_2(t)}{a_2(t)} \biggr)\vee0 \biggr]^{-1/\gamma_2},
\end{eqnarray}
for $i=1,\ldots,n$. We conclude that if we let $k = k(n)$ denote an
intermediate sequence, that is, $k \to\infty$ and $k/n \to0$ as $n
\to\infty$, then
%
%e11 #&#
\begin{equation}
\label{nine} R_n(x,y) := \frac{n}{k} P \bigl(
\widetilde{X}_i(n/k) \leq x, \widetilde {Y}_i(n/k) \leq
y \bigr) \to R(x,y)
\end{equation}
as $n \to\infty$, for all $(x,y) \in[0,\infty)^2$.

We estimate $R_n$ and hence $R$ by replacing the unknown\vspace*{1pt} quantities
$a_j(n/k)$, $b_j(n/k)$ and $\gamma_j$, $j=1,2$, by appropriate
estimators $\widehat{a}_j(n/k)$, $\widehat{b}_j(n/k)$ and $\widehat
{\gamma}_j$, and the probability $P$ by the corresponding empirical
measure. We define, therefore,
%
%e12 #&#
\begin{eqnarray}
\widehat{X}_i(n/k) &=& \biggl[ \biggl( 1 + \widehat{
\gamma}_1 \frac{X_i -
\widehat{b}_1(n/k)}{\widehat{a}_1(n/k)} \biggr)\vee0 \biggr]^{-1/\widehat
{\gamma}_1},
\nonumber
\\[-8pt]
\label{XhatYhat}
\\[-8pt]
\nonumber
\widehat{Y}_i(n/k) &=& \biggl[ \biggl( 1 + \widehat{
\gamma}_2 \frac{Y_i -
\widehat{b}_2(n/k)}{\widehat{a}_2(n/k)} \biggr)\vee0 \biggr]^{-1/\widehat
{\gamma}_2}
\end{eqnarray}
and
%
%e13 #&#
\begin{equation}
\label{Rhat} \widehat{R}_n(x,y) = \frac{1}{k} \sum
_{i=1}^n \mathbh{1}_{\{ \widehat
{X}_i(n/k) \leq x,  \widehat{Y}_i(n/k) \leq y \}}
\end{equation}
for $(x,y) \in[0,\infty)^2$; cf.~\citet{dehaanresnick1993}.

We consider the empirical process
%
%e14 #&#
\begin{equation}
\label{etan} \eta_n(x,y) = \sqrt{k} \bigl[\widehat{R}_n(x,y)
- R(x,y) \bigr], \qquad  (x,y) \in[0,\infty)^2.
\end{equation}
We will establish the asymptotic behavior of $\eta_n$ on $[\delta
,T]^2$, for any $0 < \delta< T < \infty$, but we introduce some
definitions and assumptions first. Note that from now on we will omit
the arguments $(n/k)$ where appropriate, for ease of notation.

Let $V_R(x,y)$ denote a Wiener process on $[0,\infty]^2 \setminus\{
(\infty,\infty)\}$ with ``time'' $R$, that is, a zero-mean Gaussian
process with covariance
\[
E \bigl[V_R(x,y)V_R\bigl(x',y'
\bigr) \bigr] = R\bigl(x \wedge x', y \wedge y'\bigr).
\]
Also write [cf. (\ref{nine})]
%
%e15 #&#
\begin{equation}
\label{Tn.def} T_n(x,y) = \frac{1}{k} \sum
_{i=1}^n \mathbh{1}_{\{\widetilde{X}_i \leq
x, \widetilde{Y}_i \leq y \}}, \qquad (x,y)
\in[0,\infty)^2.
\end{equation}
It is known, by \citet{einmahletal1997}, Lemma~3.1, that $\sqrt
{k}(T_n - R_n) \Rightarrow V_R$ in $D([\delta,T]^2)$, where
``$\Rightarrow$'' denotes weak convergence and $D([\delta,T]^2)$
denotes the Skorohod space of functions defined on $[\delta,T]^2$.

In order to leave the estimators $\widehat{a}_j$, $\widehat{b}_j$ and
$\widehat{\gamma}_j$, $j=1,2$, general at this stage, we simply assume
that they are chosen in such a way that:

\begin{longlist}[A3.]
\item[A1.] For some 6-variate random vector $(A_1, A_2,
B_1, B_2, \Gamma_1, \Gamma_2)$, we have the joint weak convergence
%
%e16 #&#
\begin{eqnarray}
&& \sqrt{k} \biggl( T_n - R_n,
\frac{\widehat{a}_1}{a_1}-1, \frac
{\widehat{a}_2}{a_2}-1, \frac{\widehat{b}_1 - b_1}{a_1},
\frac{\widehat
{b}_2 - b_2}{a_2}, \widehat{\gamma}_1 - \gamma_1,
\widehat{\gamma}_2 - \gamma_2 \biggr)
\nonumber
\\[-8pt]
\label{joint.weak}
\\[-8pt]
\nonumber
&& \qquad \Rightarrow (V_R,A_1,A_2,B_1,B_2,
\Gamma_1,\Gamma_2 ) %
\end{eqnarray}
in $D([\delta,T]^2) \times\mathbb{R}^6$.

Assumption A1 is fulfilled for, for example, the moment estimators of
$\gamma_j$, $a_j$ and $b_j$, provided that $k$ is chosen appropriately;
see \citet{dehaanferreira2006}, Sections~4.2 and 3.5. We further
assume the following:

\item[A2.]  The partial derivatives
\[
R_{(1)}(x,y) :=\frac{ \partial R}{\partial x}, \qquad R_{(2)}(x,y) :=
\frac{\partial R}{\partial y}
\]
exist and are continuous on $(0,\infty)^2$.

\item[A3.]  The sequence $k$ is chosen such that
\[
\sqrt{k} \sup_{(x,y) \in[\delta/2, T+1]^2} \bigl|R_n(x,y) - R(x,y)\bigr| \to0.
\]

Finally, for $j=1,2$, we define the following functions on $(0,\infty)$:
%
%e17 #&#
\begin{eqnarray}
f_j(x) &=& %
\cases{\ds \frac{x (x^{\gamma_j}-1)}{\gamma_j}, & $\quad
\gamma_j \neq0$, \vspace*{3pt}
\cr
x \log x, & $\quad
\gamma_j = 0$,}
\nonumber
\\
\label{twelve} g_j(x) &=& -x^{\gamma_j + 1},
\\
h_j(x) &=&
\cases{ \ds\frac{x (1-x^{\gamma_j})}{\gamma_j^2} +
\frac{x \log
x}{\gamma_j}, & $\quad \gamma_j \neq0$, \vspace*{3pt}
\cr
\ds -
\bigl(x \log^2 x\bigr)/2, & $\quad\gamma_j = 0$.}
\nonumber
\end{eqnarray}
We are now ready to state the basic convergence result for $\eta_n$.
\end{longlist}

%th2.1 #&#
\begin{thm}\label{etan.asymp}
Let $0 < \delta< T < \infty$. If assumptions \textup{A1}--\textup{A3} hold, then
%
%e18 #&#
%e19 #&#
\begin{eqnarray}
\eta_n(x,y) &\Rightarrow & V_R(x,y) +
R_{(1)}(x,y) \bigl[f_1(x)A_1 +
g_1(x)B_1 + h_1(x)\Gamma_1\bigr]
\nonumber
\\
\label{eleven} &&{}+ R_{(2)}(x,y) \bigl[f_2(y)A_2
+ g_2(y)B_2 + h_2(y)\Gamma_2
\bigr]
\\
& =:& \eta(x,y)\nonumber
\end{eqnarray}
in $D([\delta,T]^2)$.
\end{thm}

%re1 #&#
\begin{rem*}
Note that we take $\delta>0$, since the result does not
hold true in general for $\delta= 0$: the functions in (\ref{twelve})
are unbounded near zero for $\gamma_j < -1$. This theorem is very
similar to Theorem~5.1 in \citet{dehaanresnick1993}, where instead of
$R$ the stable tail dependence function $l(x,y)=x+y-R(x,y)$ is
estimated. We nevertheless offer a detailed proof of Theorem~\ref
{etan.asymp} in \citet{supp}, since the statement and proof of
Theorem~5.1 in \citet{dehaanresnick1993} are not completely correct;
in particular, our $\delta$ is taken to be 0 there.
\end{rem*}

%s2.1 #&#
\subsection{Parametric empirical process}\label{subsec.param}

Now suppose that the tail copula $R$ is a member of some parametric
family of tail copulas, $\mathcal{R} = \{R_{\bolds{\theta}} \dvtx
\bolds{\theta} \in\Theta\}$, where $\Theta$ is an open subset of
$\mathbb{R}^d$. Then there is a ${\bolds{\theta}}_0 = (\theta
_{01}, \ldots, \theta_{0d})^\top\in\Theta$ such that $R =
R_{\bolds{\theta}_0}$. Let $\widehat{\bolds{\theta}} =
(\widehat{\theta}_1, \ldots, \widehat{\theta}_d)^\top$ denote an
estimator of $\bolds{\theta}_0$, and consider the empirical process
%
%e20 #&#
\begin{equation}
\label{etan.hat}
\widehat{\eta}_n(x,y) = \sqrt{k} \bigl[
\widehat{R}_n(x,y) - R_{\widehat
{\bolds{\theta}}}(x,y) \bigr], \qquad (x,y) \in[0,
\infty)^2,
\end{equation}
the parametric version of (\ref{etan}). Our next result will establish
the asymptotic behavior of $\widehat{\eta}_n$. Since
%
%e21 #&#
\begin{equation}
\label{etan.decomp}
\widehat{\eta}_n(x,y) = \eta_n(x,y) +
\sqrt{k} \bigl[ R_{\bolds
{\theta}_0}(x,y) - R_{\widehat{\bolds{\theta}}}(x,y) \bigr],
\end{equation}
the asymptotic behavior of $\widehat{\eta}_n$ is an easy consequence of
Theorem~\ref{etan.asymp}, under proper assumptions. We state those
assumptions below.

\begin{longlist}[B3.]
\item[B1.]  There is a $(6+d)$-variate random vector
$(A_1,A_2,B_1,B_2,\Gamma_1,\Gamma_2,\bolds{\zeta})$ such that
%
%e22 #&#
\begin{eqnarray}
&& \sqrt{k} \biggl( T_n - R_n,
\frac{\widehat{a}_1}{a_1}-1, \frac
{\widehat{a}_2}{a_2}-1, \frac{\widehat{b}_1 - b_1}{a_1},
\frac{\widehat
{b}_2 - b_2}{a_2}, \widehat{\gamma}_1 - \gamma_1,
\widehat{\gamma}_2 - \gamma_2, \bolds{
\theta}_0 - \widehat{\bolds{\theta}} \biggr)\hspace*{-20pt}
\nonumber
\\[-8pt]
\label{B1}
\\[-8pt]
\nonumber
&& \qquad \Rightarrow (V_R,A_1,A_2,B_1,B_2,
\Gamma_1,\Gamma_2, \bolds {\zeta} )
\end{eqnarray}
in $D([\delta,T]^2) \times\mathbb{R}^{6+d}$.

\item[B2.]  The first-order partial derivatives
\begin{eqnarray*}
R_{\bolds{\theta}(1)}(x,y) &=& \frac{\partial}{\partial x} R_{\bolds{\theta}}(x,y), \qquad
R_{\bolds{\theta}(2)}(x,y) = \frac{\partial}{\partial y} R_{\bolds{\theta}}(x,y),
\\
\dt{R}_{\bolds{\theta}}(x,y) &=& \biggl( \frac{\partial}{\partial
\theta_1} R_{\bolds{\theta}}(x,y),
\ldots, \frac{\partial
}{\partial\theta_d} R_{\bolds{\theta}}(x,y) \biggr)^\top
\end{eqnarray*}
exist and are continuous for $(x,y,\bolds{\theta}) \in(0,\infty
)^2 \times B(\bolds{\theta}_0)$, for some neighborhood
$B(\bolds{\theta}_0)$ of $\bolds{\theta}_0$ in $\Theta$.

\item[B3.]  The sequence $k$ is chosen such that
%
%e23 #&#
\begin{equation}
\label{B3}
\sqrt{k} \sup_{(x,y) \in[\delta/2, T+1]^2} \bigl|R_n(x,y) -
R_{\bolds
{\theta}_0}(x,y)\bigr| \to0.
\end{equation}

Note that B3 is the same as A3; we restate it here for ease of
presentation. Also note that by virtue of B2 the second term on the
right-hand side of (\ref{etan.decomp}) is asymptotically equal in
probability to
\[
\dt{R}^\top_{\bolds{\theta}_0}(x,y)\sqrt{k}(\bolds{\theta}_0
- \widehat{\bolds{\theta}}),
\]
which, by B1, converges weakly to $\dt{R}^\top_{\bolds{\theta
}_0}(x,y)\bolds{\zeta}$. Thus we obtain the following corollary to
Theorem~\ref{etan.asymp}.
\end{longlist}

%co2.2 #&#
\begin{cor}\label{etan.hat.asymp}
Let $0 < \delta< T < \infty$. If assumptions \textup{B1}--\textup{B3} hold, then
%
%e24 #&#
%e25 #&#
\begin{eqnarray}
\widehat{\eta}_n(x,y) &\Rightarrow & V_{R_{\bolds{\theta}_0}}(x,y)
+ R_{\bolds{\theta}_0(1)}(x,y) \bigl[f_1(x)A_1 +
g_1(x)B_1 + h_1(x)\Gamma _1\bigr]
\nonumber
\\
&&{}+ R_{\bolds{\theta}_0(2)}(x,y) \bigl[f_2(y)A_2 +
g_2(y)B_2 + h_2(y)\Gamma _2\bigr]
\nonumber
\\[-8pt]
\label{asymp.2}
\\[-8pt]
\nonumber
&& {}+ \dt {R}^\top_{\bolds{\theta}_0}(x,y) \bolds{\zeta}
\\
& =:& \widehat{\eta}(x,y)
\nonumber
\end{eqnarray}
in $D([\delta,T]^2)$.
\end{cor}

%s3
%s3 #&#
\section{Transforming \texorpdfstring{$\widehat{\bolds{\eta}}$}{widehateta} into a standard
Wiener process}\label{sec.trans}

The limiting process $\widehat{\eta}$ in (\ref{asymp.2}) is of the
general form
%
%e26 #&#
\begin{equation}
\label{fourteen}
\xi(x,y) = V_R(x,y) + \sum
_{j=1}^\nu Q_j(x,y)Z_j,
\end{equation}
where $V_R$ denotes a Wiener process with time $R$, $\nu$ is a fixed
integer, $Q_1, \ldots, Q_\nu$ are deterministic functions mapping
$[\delta,T]^2$ into $\mathbb{R}$ and $Z_1, \ldots, Z_\nu$ are random
variables.

It will be more convenient to consider the set-indexed version of (\ref
{fourteen}),
%
%e27 #&#
\begin{equation}
\label{two} \xi(B) = V_R(B) + \sum_{j=1}^\nu
Q_j(B)Z_j =: V_R(B) +
\mathbf{Q}^\top (B)\mathbf{Z},
\end{equation}
where $B$ is a Borel subset of $[\delta,T]^2$, $V_R$ is a set-indexed
Wiener process with time measure $R$ and $Q_1, \ldots, Q_\nu$ are
deterministic signed measures. In the right-hand side of (\ref{two}),
$\mathbf{Q}(B)$ denotes the column vector consisting of $Q_1(B), \ldots
, Q_\nu(B)$ and $\mathbf{Z}$ denotes the column vector consisting of
$Z_1, \ldots, Z_\nu$.

We will state a general transformation result about set-indexed
processes $\xi$ of the form (\ref{two}), which we will then apply to
the process $\widehat{\eta}$ in (\ref{asymp.2}). The transformation is
a suitable extension of the ``innovation martingale transform'' first
discussed in \citeauthor{khmaladze1981} (\citeyear{khmaladze1981,khmaladze1988,khmaladze1993}) in
connection with parametric goodness-of-fit testing for univariate and
multivariate distribution functions; see, in particular, \citet{khmaladze1993}, Theorem~3.9.
A good summary of the innovation
martingale transform idea can be found in \citet{koulswordson2011};
for a variety of statistical applications we refer to \citet{mckeagueetal1995}, \citet{nikabadzestute1997},
\citet{stuteetal1998}, \citeauthor{koenkerxiao2002} (\citeyear{koenkerxiao2002,koenkerxiao2006}),
\citeauthor{khmaladzekoul2004}
(\citeyear{khmaladzekoul2004,khmaladzekoul2009}), \citet{delgadoetal2005}
and \citet{dettehetzler2009}, among others.

As in \citet{khmaladze1993}, we will call a collection of subsets $\{
A_u \dvtx  0 \le u \le1\}$ of $[\delta,T]^2$ a \emph{scanning family} over
$[\delta,T]^2$ if the following hold:
\begin{longlist}[(iii)]
\item[(i)] $\operatorname{Leb}(A_0) = 0,  \operatorname{Leb}(A_1) = (T - \delta)^2$,
\item[(ii)] $A_u \subset A_{u'}$ if $u \le u'$,
\item[(iii)] $\operatorname{Leb}(A_{u'} \setminus A_u) \to0 \mbox{ if } u'
\downarrow u$,
\end{longlist}
with Leb denoting Lebesgue measure. Note that for any $j \in\{1,\ldots
,\nu\}$ and Borel subset $B$ of $[\delta,T]^2$, the function $u \mapsto
Q_j(B \cap A_u)$ generates a signed measure on~$[0,1]$.

%th3.1 #&#
\begin{thm}\label{thm.transform}
Let $\xi$ be a set-indexed process of the form (\ref{two}). Suppose
there are functions $q_j\dvtx [\delta,T]^2 \to\mathbb{R},  1 \le j \le\nu$ that are square-integrable with respect to $R$ and that satisfy
\[
Q_j(B) = \iint_B q_j(x,y)\,\dd R(x,y),
\qquad 1 \le j \le\nu,
\]
for any Borel set $B \subset[\delta,T]^2$. Let $\{A_u\dvtx  0 \le u \le1\}
$ be a scanning family over $[\delta,T]^2$. Then the process
%
%e28 #&#
\begin{equation}
\label{one}
W_R(B) = \xi(B) - \int_0^1
\mathbf{Q}^\top(B \cap A_{\dd u}) \mathbf {I}^{-1}
\bigl(A_u^c\bigr) \iint_{A_u^c}\mathbf{q}(x,y)\,
\dd\xi(x,y)
\end{equation}
is a Wiener process with time $R$, where $\mathbf{q}(x,y)$ denotes the
column vector consisting of $q_1(x,y), \ldots, q_\nu(x,y)$, and the
matrices $\mathbf{I}(A_u^c)$ are defined by
\[
\mathbf{I}\bigl(A_u^c\bigr) = \iint_{A_u^c}
\mathbf{q}(x,y) \mathbf{q}^\top(x,y) \,\dd R(x,y), \qquad u \in[0,1)
\]
and are assumed to be invertible.
\end{thm}

Now let us return to the setup of Section~\ref{subsec.param}. We state
the following assumption.

\begin{longlist}[B4.]
\item[B4.]  For each $\bolds{\theta} \in\Theta$, the
measure $R_{\bolds{\theta}}$ can be decomposed as $R_{\bolds
{\theta}} = R^{(c)}_{\bolds{\theta}} + R^{(s)}_{\bolds{\theta
}}$, where $R^{(s)}_{\bolds{\theta}}$ satisfies
$R^{(s)}_{\bolds{\theta}}([0,\infty)^2) = 0$ and
$R^{(c)}_{\bolds{\theta}}$ is absolutely continuous with respect
to the Lebesgue measure on $(0,\infty)^2$, with a positive density
$r_{\bolds{\theta}}$ that has continuous first-order partial
derivatives with respect to $x, y, \theta_1, \ldots, \theta_d$ for all
$(x,y,\bolds{\theta}) \in(0,\infty)^2 \times B(\bolds{\theta
}_0)$, for some neighborhood $B(\bolds{\theta}_0)$ of $\bolds
{\theta}_0$ in $\Theta$.
\end{longlist}

Note that B4 allows arbitrarily large masses on the ``axes at
infinity'' $\{(x,\infty) \dvtx  x \geq0 \} \cup\{(\infty,y) \dvtx  y\geq0 \}$
for $R_{\bolds{\theta}} \in\mathcal{R}$, but excludes the case
$R_{\bolds{\theta}} \equiv R^{(s)}_{\bolds{\theta}}$, which
corresponds to (strict) tail independence.

Let us define the following functions on $[\delta,T]^2$, with $f_j,
g_j$ and $h_j$ as defined in (\ref{twelve}):
\begin{eqnarray*}
Q_1(x,y) &=& R_{\bolds{\theta}_0(1)}(x,y)f_1(x), \qquad
Q_4(x,y) = R_{\bolds{\theta}_0(2)}(x,y)f_2(y),
\\
Q_2(x,y) &=& R_{\bolds{\theta}_0(1)}(x,y)g_1(x), \qquad
Q_5(x,y) = R_{\bolds{\theta}_0(2)}(x,y)g_2(y),
\\
Q_3(x,y) &=& R_{\bolds{\theta}_0(1)}(x,y)h_1(x), \qquad
Q_6(x,y) = R_{\bolds{\theta}_0(2)}(x,y)h_2(y)
\end{eqnarray*}
and
\[
Q_{6+i}(x,y) = \frac{\partial}{\partial\theta_i} R_{\bolds{\theta
}}(x,y)
\Big|_{\bolds{\theta} = \bolds{\theta}_0}, \qquad i=1,\ldots,d.
\]
Furthermore, let $q_i$ denote the Radon--Nikodym derivatives $\dd
Q_i/\dd R_{\bolds{\theta}_0}$ for $i=1,\ldots,6+d$, or more explicitly:
\begin{eqnarray*}
q_1(x,y) &=& f_1'(x) + f_1(x)
\frac{\partial}{\partial x} \log r_{\bolds{\theta}_0}(x,y),
\\
q_2(x,y) &=& g_1'(x) + g_1(x)
\frac{\partial}{\partial x} \log r_{\bolds{\theta}_0}(x,y),
\\
q_3(x,y) &=& h_1'(x) + h_1(x)
\frac{\partial}{\partial x} \log r_{\bolds{\theta}_0}(x,y),
\\
q_4(x,y) &=& f_2'(y) + f_2(y)
\frac{\partial}{\partial y} \log r_{\bolds{\theta}_0}(x,y),
\\
q_5(x,y) &=& g_2'(y) + g_2(y)
\frac{\partial}{\partial y} \log r_{\bolds{\theta}_0}(x,y),
\\
q_6(x,y) &=& h_2'(y) + h_2(y)
\frac{\partial}{\partial y} \log r_{\bolds{\theta}_0}(x,y)
\end{eqnarray*}
and
\[
q_{6+i}(x,y) = \frac{\partial}{\partial\theta_i} \log r_{\bolds
{\theta}}(x,y)
\Big|_{\bolds{\theta} = \bolds{\theta}_0}, \qquad i=1,\ldots,d.
\]
As before, $\mathbf{q}(x,y)$ will denote the column vector consisting
of $q_1(x,y), \ldots,\break  q_{6+d}(x,y)$ for $(x,y) \in[\delta,T]^2$.

We are now ready to apply Theorem~\ref{thm.transform} to $\widehat{\eta
}$ in (\ref{asymp.2}). Instead of arbitrary Borel sets $B$, we consider
rectangles $[\delta,x] \times[\delta,y] \subset[\delta, T]^2$, with
\[
\widehat{\eta}\bigl([\delta,x] \times[\delta,y]\bigr) := \widehat{\eta}(x,y)-
\widehat{\eta}(\delta,y) - \widehat{\eta}(x,\delta) + \widehat{\eta }(\delta,
\delta).
\]
We also introduce the scanning family $A_u = [\delta,T] \times[\delta
,(1-u)\delta+ uT]$ for $0 \le u \le1$ and define the corresponding matrices
%
%e29 #&#
\begin{equation}
\label{I.mat}
\mathbf{I}(t) = \int_\delta^T\!\! \int
_t^T \mathbf{q}\bigl(s',t'
\bigr) \mathbf{q}^\top \bigl(s',t'\bigr) \,
\dd R_{\bolds{\theta}_0}\bigl(s',t'\bigr), \qquad t \in[
\delta,T).
\end{equation}

\begin{rem*}
From a likelihood theory point of view, the functions
$q_1, \ldots,\break  q_{6+d}$ can be seen as score functions corresponding to
the estimated values $a_1, a_2$, $b_1, b_2, \gamma_1, \gamma_2, \theta
_{01}, \ldots, \theta_{0d}$, and the matrix $\mathbf{I}(t)$ can be seen
as a partial Fisher information matrix constructed from these score functions.
\end{rem*}

%co3.2 #&#
\begin{cor}\label{cor.transform}
If assumptions \textup{B2} and \textup{B4}, restricted to $\bolds{\theta} =
\bolds{\theta}_0$, hold, and the matrices $\mathbf{I}(t)$ in (\ref
{I.mat}) are invertible, then the process
\begin{eqnarray*}
&& W_R\bigl([\delta,x] \times[\delta,y]\bigr)
\\
&& \qquad = \widehat{\eta}\bigl([\delta,x] \times[\delta,y]\bigr)
\\
&&\quad\qquad{}- \int_\delta^x\!\! \int
_\delta^y \mathbf{q}^\top(s,t) \biggl(
\mathbf{I}^{-1}(t) \int_\delta^T \!\!\int
_t^T \mathbf{q}\bigl(s',t'
\bigr)\,\dd \widehat{\eta}\bigl(s',t'\bigr) \biggr) \,
\dd R_{\bolds{\theta}_0}(s,t)
\end{eqnarray*}
is a Wiener process with time $R_{\bolds{\theta}_0}$ on $[\delta
,T] \times[\delta,T)$.
\end{cor}

%can be uniquely defined even when $\mathbf{I}(t)$ is not of full rank.
%We refer to \citet{{tsigroshvili1998} for details.

In order to obtain a \emph{standard} Wiener process from $\widehat{\eta
}$, we normalize $W_R$ in the usual way, as follows.
%
%co3.3 #&#
\begin{cor}\label{cor.norm}
If assumptions \textup{B2} and \textup{B4}, restricted to $\bolds{\theta} =
\bolds{\theta}_0$, hold, and the matrices $\mathbf{I}(t)$ in (\ref
{I.mat}) are invertible, then the process
%
%e30 #&#
%e31 #&#
\begin{eqnarray}
&& W\bigl([\delta,x] \times[\delta,y]\bigr)
\nonumber\\
&& \qquad = \int_\delta^x\!\! \int
_\delta^y \frac{1}{\sqrt{r_{\bolds{\theta
}_0}(s,t)}} \, \dd
W_R\bigl([\delta,s] \times[\delta,t]\bigr)
\nonumber
\\[-8pt]
\label{W}
\\[-8pt]
\nonumber
&&\qquad = \int_\delta^x\!\! \int_\delta^y
\frac{1}{\sqrt{r_{\bolds{\theta
}_0}(s,t)}} \, \dd\widehat{\eta}(s,t)
\\
&&\quad\qquad{}- \int_\delta^x \!\!\int
_\delta^y \mathbf{q}^\top(s,t) \biggl(
\mathbf{I}^{-1}(t) \int_\delta^T \!\!\int
_t^T \mathbf{q}\bigl(s',t'
\bigr)\,\dd\widehat{\eta}\bigl(s',t'\bigr) \biggr)
\sqrt {r_{\bolds{\theta}_0}(s,t)} \, \dd t \, \dd s
\nonumber
\end{eqnarray}
is a standard Wiener process on $[\delta,T] \times[\delta, T)$.
\end{cor}

%s4
%s4 #&#
\section{Goodness-of-fit testing}\label{sec.gof}

In Section~\ref{sec.estimator} we introduced the parametric empirical
process $\widehat{\eta}_n$ as the normalized difference between
$R_{\widehat{\bolds{\theta}}}$ and the semi-parametric estimator
$\widehat{R}_n$, and derived its weak limit $\widehat{\eta}$. In
Section~\ref{sec.trans} we described a transformation from $\widehat
{\eta}$ into a standard Wiener process $W$. In this section, we will
apply the empirical version of the same transformation to $\widehat{\eta
}_n$, and prove that the resulting empirical process converges weakly
to a standard Wiener process. This is the main result of this paper.

Define the empirical version of $W$ in (\ref{W}) as follows, for $(x,y)
\in[\delta,T] \times[\delta,T)$:
%
%e32 #&#
\begin{eqnarray}
&& W_n\bigl([\delta,x] \times[\delta,y]\bigr)
\nonumber
\\
\label{Wn}
&&\qquad= \int_\delta^x\!\! \int_\delta^y
\frac{1}{\sqrt{r_{\widehat{\bolds
{\theta}}}(s,t)}} \, \dd\widehat{\eta}_n(s,t)
\\
&&\qquad\quad{} - \int_\delta^x \!\!\int
_\delta^y \widehat{\mathbf{q}}^\top (s,t)
\biggl( \widehat{\mathbf{I}}^{-1}(t) \int_\delta^T\!\!
\int_t^T \widehat{\mathbf{q}}
\bigl(s',t'\bigr)\,\dd\widehat{\eta}_n
\bigl(s',t'\bigr) \biggr) \sqrt{r_{\widehat{\bolds{\theta}}}(s,t)}
\,\dd t\,\dd s.
\nonumber
\end{eqnarray}
Here, the vectors $\widehat{\mathbf{q}}$ and the matrices $\widehat
{\mathbf{I}}$ are obtained by replacing the unknown marginal tail\vspace*{1pt}
indices $\gamma_1, \gamma_2$ and the unknown parameter $\bolds
{\theta}_0$ in the definition of $\mathbf{q}$ by their estimators
$\widehat{\gamma}_1, \widehat{\gamma}_2, \widehat{\bolds{\theta}}$.

For functions $\varphi\dvtx [\delta,T]^2 \to\mathbb{R}$, we introduce the seminorm
%
%e33 #&#
\begin{eqnarray}
\| \varphi\|_{\mathrm{HK}} &:=& V^{(2)} (\varphi) +
V^{(1)} \bigl(\varphi( \,\cdot\,, \delta)\bigr) + V^{(1)}
\bigl(\varphi(\delta,\,\cdot\,)\bigr)
\nonumber
\\[-8pt]
\label{hk.norm}
\\[-8pt]
\nonumber
&&{} + V^{(1)} \bigl(\varphi( \,\cdot\,, T)\bigr) +
V^{(1)} \bigl(\varphi(T,\,\cdot\,)\bigr),
\end{eqnarray}
where $V^{(1)}$ denotes the univariate total variation over $[\delta
,T]$, and $V^{(2)}$ denotes the bivariate (Vitali) total variation over
$[\delta,T]^2$, as defined in \citet{owen2005}, for example. The
seminorm $\| \cdot\|_{\mathrm{HK}}$ is sometimes called the \emph
{Hardy--Krause variation} in the literature, in recognition of \citet{hardy1905} and \citet{krause1903}.

For notational convenience, let us also denote
\begin{eqnarray*}
\rho_1(x,y,\bolds{\theta}) &=&  \frac{\partial}{\partial x} \log
r_{\bolds{\theta}}(x,y), \qquad \rho_2(x,y,\bolds{\theta}) =
\frac{\partial}{\partial y} \log r_{\bolds{\theta}}(x,y),
\\
\rho_{2+i}(x,y,\bolds{\theta}) &=& \frac{\partial}{\partial\theta
_i} \log
r_{\bolds{\theta}}(x,y), \qquad  i=1,\ldots,d
\end{eqnarray*}
and
\[
\Delta\rho_j (x,y) = \rho_j(x,y,\widehat{\bolds{
\theta}}) - \rho _j(x,y,\bolds{\theta}_0), \qquad j=1,
\ldots,2+d.
\]
Similarly, let
\[
\sigma(x,y,\bolds{\theta}) = r_{\bolds{\theta}}(x,y)^{-1/2}, \qquad \Delta
\sigma(x,y) = \sigma(x,y,\widehat{\bolds{\theta}}) - \sigma(x,y,\bolds{
\theta}_0).
\]
We introduce the following assumption:
\begin{longlist}[B5.]
\item[B5.] For $j=1,\ldots,2+d$, $\|\rho_j(x,y,\bolds
{\theta}_0)\|_{\mathrm{HK}} < \infty$ and
$\|\Delta\rho_j(x,y)\|_{\mathrm
{HK}} = o_P(1)$. Furthermore, $\|\sigma(x,y,\bolds{\theta}_0)\|
_{\mathrm{HK}} < \infty$ and $\|\Delta\sigma(x,y)\|_{\mathrm{HK}} =
o_P(1)$.
\end{longlist}

Given the consistency of $\widehat{\bolds{\theta}}$, which is
implied by B1, a sufficient (but not necessary) condition for B5 is the
existence and continuity of the partial derivatives
\[
\frac{\partial\varphi(x,\delta,\bolds{\theta})}{\partial x},   \frac{\partial\varphi(x,T,\bolds{\theta})}{\partial x},   \frac
{\partial\varphi(\delta,y,\bolds{\theta})}{\partial y},
\frac
{\partial\varphi(T,y,\bolds{\theta})}{\partial y},   \frac
{\partial^2 \varphi(x,y,\bolds{\theta})}{\partial x \,\partial y}
\]
on $(x,y, \bolds{\theta}) \in[\delta,T]^2 \times B(\bolds
{\theta}_0)$, for some neighborhood $B(\bolds{\theta}_0)$ of
$\bolds{\theta}_0$ in $\Theta$, for $\varphi= \sigma$ and
$\varphi= \rho_j$, $j=1,\ldots,2+d$.

We can now present the main result of this paper. %Its proof is given
%in Section~\ref{sec.proofs}.
%
%th4.1 #&#
\begin{thm}\label{thm.main}
Let $0 < \delta< \tau< T$, and let $W$ and $W_n$ be defined as in
(\ref{W}) and (\ref{Wn}). If assumptions \textup{B1}--\textup{B5} hold, then
\[
W_n\bigl([\delta,x] \times[\delta,y]\bigr) \Rightarrow W\bigl([
\delta,x] \times [\delta,y]\bigr)
\]
in $D([\delta,\tau]^2)$.
\end{thm}

Note that Theorem~\ref{thm.main} yields that under the null hypothesis
$R \in\mathcal{R}$, we obtain a \emph{distribution-free} limiting
process $W$ (a standard bivariate Wiener process). Hence $W_n$ can be
used as a ``test process'' for producing a myriad of asymptotically
distribution-free test statistics to test this null hypothesis. We
will consider examples of such tests in Section~\ref{sec.sim}.

\begin{rem*}
By taking $\mathcal{R} = \{R_0\}$, where $R_0$ is a
fully specified tail copula, we can use Theorem~\ref{thm.main} for
testing the null hypothesis $R = R_0$. In this case, the process
$\widehat{\eta}_n$ in the definition of $W_n$ [see (\ref{Wn})]\vspace*{1pt} reduces
to $\eta_n$ as defined in $(\ref{etan})$, $r_{\widehat{\bolds
{\theta}}}$~reduces to $r_0 = \dd R_0/ \dd\operatorname{Leb}$ and $\widehat
{\mathbf{q}}$ and $\widehat{\mathbf{I}}$ are determined by $r_0$ and
$\widehat{\gamma}_1, \widehat{\gamma}_2$. We will consider an example
of testing $R=R_0$ in Section~\ref{sec.sim}.
\end{rem*}

%s5
%s5 #&#
\section{Multivariate extension}\label{sec.multi}

In this section we extend Theorem~\ref{thm.main} from the bivariate to
the $m$-dimensional setting, for $m>2$. The proof will be omitted, but
it follows very similar lines as in the bivariate case. In particular,
Theorem~\ref{thm.transform} immediately generalizes to dimension $m$
and then serves as a basis for the main result of this section.
%We will consider the multivariate case also in the simulation section.

So suppose that we have an i.i.d. sample $\mathbf{X}_1, \ldots, \mathbf
{X}_n$ from an $m$-variate d.f. $F$ with marginal d.f.'s $F_1, \ldots, F_m$.
We write, for each $i \in\{1,\ldots,n\}$, $\mathbf{X}_i = (X_{i1},
\ldots, X_{im})^\top$, where $X_{ij}$ has d.f. $F_j$. We assume that $F$
is in the max-domain of attraction of an $m$-variate extreme value d.f.
$G$, so there exist normalizing sequences $a_1(n), \ldots, a_m(n) >0$
and $b_1(n), \ldots, b_m(n) \in\mathbb{R}$ such that
\[
P \biggl( \frac{\bigvee_{i=1}^n X_{i1} - b_1(n)}{a_1(n)} \le x_1, \ldots ,
\frac{\bigvee_{i=1}^n X_{im} - b_m(n)}{a_m(n)} \le x_m \biggr) \stackrel{d} {\rightarrow} G(
\mathbf{x}),
\]
with $\mathbf{x} = (x_1, \ldots, x_m)^\top\in\mathbb{R}^m$. We
assume, as in the bivariate case, that the sequences $a_j$ and $b_j$,
$j=1,\ldots,m$, are chosen in such a way that $G$ has marginal d.f.'s of
the form
\[
G_j(x) = \exp\bigl\{-(1+\gamma_jx)^{-1/\gamma_j}
\bigr\}, \qquad 1+\gamma_jx > 0,
\]
for some $\gamma_1, \ldots, \gamma_m \in\mathbb{R}$. We will denote
$\bolds{\gamma} = (\gamma_1, \ldots, \gamma_m)^\top$. The d.f. $G$
is then characterized by the marginal tail indices $\bolds{\gamma
}$ and the $m$-variate tail copula
\[
R(\mathbf{x}) := \lim_{t \to\infty} t P \Biggl( \bigcap
_{j=1}^m \bigl\{ 1 - F_j(X_{1j})
\le x_j/t \bigr\} \Biggr), \qquad \mathbf{x} \in[0,
\infty]^m \setminus\{\bolds{\infty}\},
\]
where $\bolds{\infty}$ denotes the point $(\infty, \ldots, \infty
)$.

\begin{rem*}
In the remainder of this section we consider $R$ defined
on the restricted domain $[0,\infty)^m$ [cf. (\ref{R.def})] because
our processes and transformations are not defined outside this region.
The bivariate tail copula $R$ defined on $[0,\infty)^2$ determines $R$
on the full domain $[0,\infty]^2 \setminus\{(\infty, \infty)\}$. In
contrast, for $m>2$ the tail copula $R$ defined on $[0,\infty)^m$ in
general does \emph{not} determine $R$ on the full domain $[0,\infty]^m
\setminus\{\bolds{\infty}\}$.
\end{rem*}

Let $\mathcal{R} = \{R_{\bolds{\theta}}\dvtx  \bolds{\theta} \in
\Theta\}$ denote a parametric family of $m$-variate tail copulas on
$[0,\infty)^m$, parametrized by $\bolds{\theta} = (\theta_1, \ldots
, \theta_d)^\top\in\Theta$, an open subset of $\mathbb{R}^d$. Our aim
is to enable the construction of tests for the null hypothesis $R \in
\mathcal{R}$ against the alternative $R \notin\mathcal{R}$.

For fixed $\bolds{\theta} \in\Theta$, $R_{\bolds{\theta}}$
can be seen as an equivalence class of tail dependence structures
(i.e., tail copulas defined on the full domain) containing one or more
elements. Under the additional assumption that $R_{\bolds{\theta
}}$ puts no mass on $[0,\infty]^m \setminus(\{\bolds{\infty}\}
\cup[0,\infty)^m) $, $R_{\bolds{\theta}}$ contains exactly one
element (as in the bivariate case).

Suppose the null hypothesis holds true, with $R = R_{\bolds{\theta
}_0}$, for some $\bolds{\theta}_0 \in\Theta$. Let $\widehat
{\bolds{\theta}}$ denote an estimator for $\bolds{\theta}_0$.
As in Section~\ref{sec.estimator}, we let $k = k(n)$ denote an
intermediate sequence and define the parametric empirical process
\[
\widehat{\eta}_n(\mathbf{x}) = \sqrt{k} \bigl[\widehat{R}_n(
\mathbf{x}) - R_{\widehat{\bolds{\theta}}}(\mathbf{x}) \bigr], \qquad \mathbf{x} \in[0,
\infty)^m,
\]
where
\[
\widehat{R}_n(\mathbf{x}) = \frac{1}{k} \sum
_{i=1}^n \mathbh{1}_{\bigcap
_{j=1}^m \{ \widehat{X}_{ij}(n/k) \le x_j\}}, \qquad \mathbf{x}
\in [0,\infty)^m,
\]
with $\widehat{X}_{ij}(n/k), (i,j) \in\{1,\ldots,n\} \times\{1,\ldots
,m\}$, defined similarly as in (\ref{XhatYhat}). Let $R_n$ and $T_n$
denote the obvious $m$-variate extensions of (\ref{nine}) and (\ref
{Tn.def}), let $0 < \delta< T < \infty$ and let C1--C4 denote the
natural $m$-variate extensions of assumptions B1--B4 of Sections~\ref{sec.estimator} and \ref{sec.trans}.

To state the analog of assumption B5 for the $m$-variate case, we
extend the seminorm (\ref{hk.norm}) to $m$-variate functions by
induction, as follows: For any function $\varphi\dvtx [\delta,T]^m \to
\mathbb{R}$, and $i \in\{1, \ldots, m\}$, we define $\varphi_{\delta
,i}\dvtx [\delta,T]^{m-1} \to\mathbb{R}$ to be the restriction of $\varphi$
to the subset of $[\delta,T]^m$ with the $i$th coordinate fixed at
$\delta$, and we define $\varphi_{T,i}$ analogously. Then we let
%
%e34 #&#
\begin{equation}
\label{hk.norm.multi}
\| \varphi\|_{\mathrm{HK}}^{(m)} := V^{(m)} (
\varphi) + \sum_{i=1}^m \|
\varphi_{\delta,i} \|_{\mathrm{HK}}^{(m-1)} + \sum
_{i=1}^m \| \varphi _{T,i}
\|_{\mathrm{HK}}^{(m-1)},
\end{equation}
with $V^{(m)}$ denoting the $m$-variate (Vitali) total variation over
$[\delta,T]^m$ and $\|\varphi\|_{\mathrm{HK}}^{(2)}$ as defined in (\ref
{hk.norm}). We also let $\rho_j, \Delta\rho_j, \sigma, \Delta\sigma$
be defined as in Section~\ref{sec.gof}, for $j=1,\ldots,m+d$.

\begin{longlist}[C5.]
\item[C5.]  For $j=1,\ldots,m+d$, $\|\rho_j(\mathbf
{x},\bolds{\theta}_0)\|_{\mathrm{HK}}^{(m)} < \infty$ and $\|\Delta
\rho_j(\mathbf{x})\|_{\mathrm{HK}}^{(m)} = o_P(1)$. Furthermore, $\|\sigma
(\mathbf{x},\bolds{\theta}_0)\|_{\mathrm{HK}}^{(m)} < \infty$ and $\|
\Delta\sigma(\mathbf{x})\|_{\mathrm{HK}}^{(m)} = o_P(1)$.
\end{longlist}

Now, let us introduce the functions $Q_j$ and $q_j = \dd Q_j/\dd
R_{\bolds{\theta}_0}$, for $j=1,\ldots, 3m+d$, as the natural
$m$-variate extensions of the bivariate functions introduced before
Corollary~\ref{cor.transform}, and let us denote by $\mathbf{q}(\mathbf
{x})$ the column vector consisting of $q_1(\mathbf{x}), \ldots,q_{3m+d}(\mathbf{x})$. Further, let us write $[\bolds{\delta},
\mathbf{x}] = [\delta, x_1] \times\cdots\times[\delta, x_m]$, $S_t =
[\delta,T]^{m-1} \times(t,T]$, and introduce matrices
\[
\mathbf{I}(t) = \int_{S_t} \mathbf{q}(\mathbf{s})
\mathbf{q}^\top (\mathbf{s}) \,\dd R_{\bolds{\theta}_0}(\mathbf{s}), \qquad t
\in [\delta,T),
\]
which are assumed to be invertible. Then the $m$-variate analog of the
transformed empirical process $W_n$ in (\ref{Wn}) is
%
%e35 #&#
%e36 #&#
\begin{eqnarray}
W_n\bigl([\bolds{\delta}, \mathbf{x}]\bigr)
&=& \int_{[\bolds{\delta},
\mathbf{x}]} \frac{1}{\sqrt{r_{\widehat{\bolds{\theta}}}(\mathbf
{s})}} \, \dd\widehat{
\eta}_n(\mathbf{s})
\nonumber
\\[-8pt]
\label{Wn.multi}
\\[-8pt]
\nonumber
&&{} - \int_{[\bolds{\delta}, \mathbf{x}]} \widehat {
\mathbf{q}}^\top(\mathbf{s}) \biggl( \widehat{\mathbf{I}}^{-1}(t)
\int_{S_t} \widehat{\mathbf{q}}\bigl(\mathbf{s}'
\bigr)\,\dd\widehat{\eta}_n\bigl(\mathbf{s}'\bigr)
\biggr) \sqrt{r_{\widehat{\bolds{\theta}}}(\mathbf{s})} \, \dd \mathbf{s},
\end{eqnarray}
where $\widehat{\mathbf{q}}$ and $\widehat{\mathbf{I}}$ are obtained by
replacing $\bolds{\gamma}$ and $\bolds{\theta}_0$ by $\widehat
{\bolds{\gamma}}$ and $\widehat{\bolds{\theta}}$ in the
definition of~$\mathbf{q}$.

We are now ready to state the multivariate analog of Theorem~\ref
{thm.main}. As in the bivariate case, this result can be used as a
basis for producing a multitude of asymptotically distribution-free
goodness-of-fit tests for a parametric model $\mathcal{R}$ (as well as
for a fully specified tail copula $R_0$).
%
%th5.1 #&#
\begin{thm}\label{theorem4.multi}
Let $m>2$. Furthermore, let $0 < \delta< \tau< T$, and let $W_n$ be
defined as in (\ref{Wn.multi}). If assumptions \textup{C1}--\textup{C5} hold, then
\[
W_n\bigl([\bolds{\delta},\mathbf{x}]\bigr) \Rightarrow W\bigl([
\bolds{\delta },\mathbf{x}]\bigr)
\]
in $D([\delta,\tau]^m)$, where $W$ is a standard $m$-variate Wiener process.
\end{thm}

%s6
%s6 #&#
\section{Simulation study}\label{sec.sim}

In this section we consider some specific functionals of $W_n$ under
the null and alternative hypotheses, for three bivariate models
$\mathcal{R}$. We will see in Monte Carlo simulations that under the
null hypothesis our limit theorems yield good approximations for finite
sample size $n$, and we also find that the resulting tests have good
power properties. This shows the applicability of our method.

The three models we consider are the following:

\begin{longlist}[Model 2.]
\item[\textit{Model} 1.]  $R(x,y) = x + y - \sqrt{x^2 + y^2}$;

\item[\textit{Model} 2.] $R \in\mathcal{R} =  \{R_\theta\dvtx  R_\theta(x,y) = x
+ y -  (x^{1/\theta} + y^{1/\theta} )^\theta, \theta\in(0,1)
\}$;

\item[\textit{Model} 3.]  $R \in\mathcal{R} =  \{R_\psi\dvtx  R_\psi(x,y) = \psi
 (x + y - \sqrt{x^2 + y^2} ), \psi\in(0,1) \}$.
 \end{longlist}

Model 2 is the widely used \emph{logistic family} of tail copulas.
Model 1 is a fully specified tail copula and a special case of Model 2.
Model 3 is a mixture between Model 1 and the tail independence model
($R \equiv0$). Note that the tail copulas of Model 3 assign mass to
the axes at infinity; indeed, the parameter $\psi$ determines how much
mass is assigned there.

For each model, we first generate 300 samples of size $n=1500$ from a
``null hypothesis d.f.'' $F_0$ for which the model is correct. We use
these samples to assess the finite-sample performance of our main
convergence result, Theorem~\ref{thm.main}. Next, we generate, for each
model, 100 samples of size $n=1500$ from an ``alternative hypothesis
d.f.'' $F_a$ for which the model is incorrect. These samples are used for
power calculations.

In Section~\ref{sec.sim.1} below, we present the data generating
distributions used for each model. Then in Section~\ref{sec.sim.2}
we describe our simulation results.
%how these samples are used to assess the quality of the approximation
%provided by Theorem~\ref{thm.main} and to carry out power calculations.
Additional details about the simulations, including the verification of
assumptions and the computer code that was used, can be found in \citet{supp}.

%s6.1 #&#
\subsection{Data generating distributions}\label{sec.sim.1}
To test for Models 1 and 2 under the null hypothesis, we generate
samples from the bivariate Cauchy distribution on the positive quadrant
with density
%
%e37 #&#
\begin{equation}
\label{biv.cauchy}
f_0(x,y) = \frac{2}{\pi(1+x^2+y^2)^{3/2}}, \qquad (x,y) \in[0,
\infty)^2.
\end{equation}
This distribution satisfies Model 1, and therefore also Model 2, with
$\theta= 1/2$.

To test for Model 3 under the null hypothesis, we sample from the
bivariate mixture random vector
%
%e38 #&#
\begin{equation}
\label{mixture}
\bigl(IX_1 + (1-I)X_2, IY_1 +
(1-I)Y_2\bigr),
\end{equation}
where $I, (X_1,Y_1), (X_2,Y_2)$ are independent, $I \sim\operatorname{Bernoulli}(0.75)$, $(X_1, Y_1)$ has the bivariate Cauchy distribution
(\ref{biv.cauchy}) on the positive quadrant and $(X_2, Y_2)$ is a pair
of standard Cauchy absolute values coupled by the counter\-monotonic
copula. Since $(X_1, Y_1)$ has the Model 1 tail copula and $(X_2, Y_2)$
has tail independence, mixture~(\ref{mixture}) has the Model 3 tail
copula with $\psi=0.75$.

To test for Model 1 under the alternative hypothesis, we sample from a
mixture random vector as in (\ref{mixture}), where $I, (X_1,Y_1),
(X_2,Y_2)$ are independent and $I \sim\operatorname{Bernoulli}(0.75)$ as
before, but $(X_1, Y_1)$ has a bivariate logistic d.f. with Fr\'{e}chet marginals,
%
%e39 #&#
\begin{equation}
\label{biv.logist}
\hspace*{13pt} F(x,y) = \exp \bigl\{- \bigl[(1+x)^{-4} +
(1+y)^{-4} \bigr]^{1/4} \bigr\}, \qquad (x,y) \in(-1,
\infty)^2,
\end{equation}
and $(X_2,Y_2)$ has identical marginal d.f.'s as in (\ref{biv.logist}),
coupled by the countermonotonic copula. The resulting d.f. has the tail copula
\[
R(x,y) = 0.75 \bigl[x+y-\bigl(x^4 + y^4
\bigr)^{1/4} \bigr], \qquad (x,y) \in[0,\infty)^2.
\]

To test for Model 2 under the alternative hypothesis, we sample from
the bivariate vector
%
%e40 #&#
\begin{equation}
\label{factor}
\bigl(\lambda Z_1 + (1-\lambda) Z_2, \mu
Z_1 + (1-\mu) Z_2\bigr),
\end{equation}
where $Z_1$ and $Z_2$ denote independent standard Pareto random
variables, and $\lambda, \mu\in(0,1)$ are deterministic coefficients.
We set $\lambda= 0.95$, $\mu= 0.65$ for the simulations. The random
vector (\ref{factor}) is a simple example of the linear factor model,
with associated tail copula
\[
R(x,y) = \min\{\lambda x, \mu y\} + \min\bigl\{(1-\lambda)x, (1-\mu) y\bigr\},
\qquad (x,y) \in[0,\infty)^2.
\]

Finally, to test for Model 3 under the alternative hypothesis, we
sample from the following asymmetric logistic d.f. with Fr\'{e}chet marginals:
%
%e41 #&#
%e42 #&#
\begin{eqnarray}
F_a(x,y) &=&  \exp \biggl\{- \biggl[
\frac{1-\phi}{1+y} + \sqrt{\frac
{1}{(1+x)^2} + \frac{\phi^2}{(1+y)^2}} \biggr] \biggr
\},
\nonumber
\\[-8pt]
\label{asym.logist}
\\[-8pt]
\eqntext{(x,y) \in(-1,\infty)^2,}
\end{eqnarray}
with $\phi= 0.25$. This d.f. has the tail copula
\[
R(x,y) = x + \phi y - \sqrt{x^2 + (\phi y)^2},
\qquad (x,y) \in[0,\infty)^2.
\]

%s6.2 #&#
\subsection{Simulation results}\label{sec.sim.2}
From each generated sample, the empirical process $W_n([\delta,x]
\times[\delta,y])$ of (\ref{Wn}) is computed on a $200 \times200$
grid $\mathcal{G}$ of uniform mesh length spanned over $[\delta,\tau
]^2$, with $\delta=0.001$ and $\tau=1.001$. We take $k=250$ and $T=2$
for all computations. The estimators $\widehat{\gamma}_j$ and $\widehat
{a}_j$, $j=1,2$, are taken to be the moment estimators [see, e.g.,
\citet{dehaanferreira2006}, Sections~4.2~and~3.5], and we set as
usual $\widehat{b}_1 = X_{n-k: n}$, $\widehat{b}_2 = Y_{n-k: n}$, with
$X_{i:n}$, $Y_{i:n}$ denoting the marginal order statistics. To
estimate the parameters $\theta$ and $\psi$ of Models 2 and 3, we use
the method of moments estimator described in \citet{einmahletal2008},
with auxiliary function $g \equiv1$.

To compare the process $W_n$ to a standard Wiener process, three test
statistics are computed from each path of $W_n$. These are:
\begin{eqnarray*}
\kappa_n &=& \max_{(x,y) \in\mathcal{G}} \bigl|W_n\bigl([
\delta,x] \times [\delta,y]\bigr) \bigr| \qquad \mbox{(Kolmogorov--Smirnov type)},
\\
\omega_n^2 &=& \|\mathcal{G}\|^2 \sum
_{(x,y) \in\mathcal{G}} W_n\bigl([\delta ,x] \times[\delta,y]
\bigr)^2 \qquad \mbox{(Cram\'{e}r--von Mises type)},
\\
A_n^2 &=& \|\mathcal{G}\|^2\sum
_{(x,y) \in\mathcal{G}} \frac{W_n([\delta
,x] \times[\delta,y])^2}{(x-\delta)(y-\delta)}\qquad \mbox{(Anderson--Darling type),}
\end{eqnarray*}
where $\|\mathcal{G}\|$ denotes the mesh length of the grid $\mathcal
{G}$, that is, $1/200$. To create benchmark distribution tables for these
statistics, we also simulate 10{,}000 true standard Wiener process paths
on the grid $\mathcal{G}$, and we compute the same test statistics for
each path. We denote these statistics, computed from the true standard
Wiener process, by $\kappa$, $\omega^2$ and $A^2$. In view of the
asymptotically distribution-free nature of our approach, these
benchmark tables need to be produced only once.
%f1 #&#
\begin{figure}[t]

\includegraphics{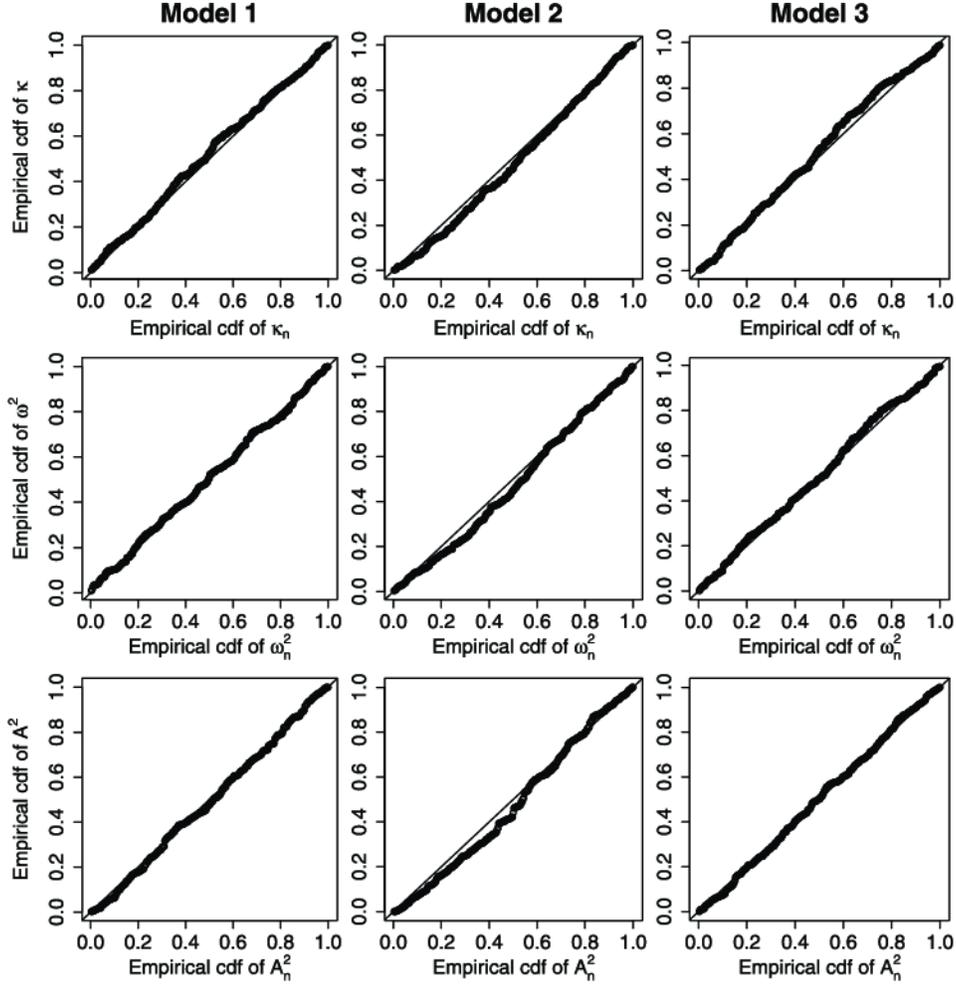}

\caption{PP-plots for the Kolmogorov--Smirnov, Cram\'{e}r--von Mises
and Anderson--Darling type test statistics.}
\label{fig.pp}
\end{figure}

For the 300 values of $\kappa_n$, $\omega^2_n$ and $A^2_n$ computed
from the null hypothesis samples, we construct PP-plots to compare
their empirical d.f.'s with the empirical d.f.'s of $\kappa$, $\omega^2$ and
$A^2$, respectively. The results are shown in Figure~\ref{fig.pp}. We
see a good match of empirical d.f.'s for all three models, which shows
that Theorem~\ref{thm.main} yields good finite-sample approximations.
This is also confirmed by the empirical size table given in the left
panel of Table~\ref{table.size}, where the observed fractions of Type I errors at the
5\% significance level are shown. Note that these numbers are
consistent with draws from a Binomial$(300,0.05)$ distribution.

For the 100 values of the test statistics computed under each
alternative hypothesis, we present the observed fraction of rejections
at the 5\% significance level in the right panel of Table~\ref{table.size}. All three
tests have quite high power.

%
%t1 #&#
\begin{table}
\caption{Observed rejection frequencies at the 5\% significance level
under null and alternative hypotheses}
\label{table.size}
\begin{tabular*}{\tablewidth}{@{\extracolsep{\fill}}lcccccc@{}}
\hline
& \multicolumn{3}{c}{\textbf{Null}} & \multicolumn{3}{c@{}}{\textbf{Alternative}}\\[-4pt]
& \multicolumn{3}{c}{\hrulefill} & \multicolumn{3}{c@{}}{\hrulefill}
\\
& \textbf{Model 1} & \textbf{Model 2} & \textbf{Model 3}
& \textbf{Model 1} & \textbf{Model 2} & \textbf{Model 3}\\
\hline
$\kappa_n$   & $15/300$ & $19/300$ & \phantom{0}$9/300$  & \phantom{0}$97/100$  & $92/100$ & \phantom{0}$97/100$ \\[2pt]
$\omega_n^2$ & $16/300$ & $11/300$ & $13/300$ & \phantom{0}$99/100$  & $90/100$ & \phantom{0}$97/100$ \\[2pt]
$A_n^2$      & $21/300$ & $17/300$ & $18/300$ & $100/100$ & $95/100$ & $100/100$
\\
\hline
\end{tabular*}
\end{table}

%s7
%s7 #&#
\section{Proofs}\label{sec.proofs}
\mbox{}
\begin{pf*}{Proof of Theorem~\protect\ref{thm.transform}}
First note that the
terms following $V_R(B)$ in (\ref{two}) are ``annihilated'' by the
transformation (\ref{one}):
\begin{eqnarray*}
&& \mathbf{Q}^\top(B)\mathbf{Z} - \int_0^1
\mathbf{Q}^\top(B \cap A_{\dd u}) \mathbf{I}^{-1}
\bigl(A_u^c\bigr)\iint_{A_u^c}\mathbf{q}(x,y)
\mathbf {q}^\top(x,y)\,\dd R(x,y)\mathbf{Z}
\\
&&\qquad = \mathbf{Q}^\top(B)\mathbf{Z} - \int_0^1
\mathbf{Q}^\top(B \cap A_{\dd u}) \mathbf{Z} = 0.
\end{eqnarray*}
Thus we can now compute, for Borel sets $B, B' \subset[\delta,T]^2$,
\begin{eqnarray*}
&& \cov \bigl[W_R(B)W_R\bigl(B'
\bigr) \bigr]
\\
&& \qquad = E \biggl[ V_R(B) - \int_0^1
\mathbf{Q}^\top(B \cap A_{\dd u}) \mathbf {I}^{-1}
\bigl(A_u^c\bigr)\iint_{A_u^c}\mathbf{q}(x,y)\,\dd
V_R(x,y) \biggr]
\\
&& \qquad\quad{}\times \biggl[ V_R\bigl(B'\bigr) - \int
_0^1 \mathbf{Q}^\top
\bigl(B' \cap A_{\dd u'}\bigr) \mathbf{I}^{-1}
\bigl(A_{u'}^c\bigr)\iint_{A_{u'}^c}\mathbf{q}(x,y)\,
\dd V_R(x,y) \biggr]
\\
&&\qquad = R\bigl(B \cap B'\bigr)
\\
&&\quad \qquad{} - \int_0^1 \mathbf{Q}^\top(B
\cap A_{\dd u}) \mathbf {I}^{-1}\bigl(A_u^c
\bigr) \mathbf{Q}\bigl(B' \cap A_u^c\bigr)
\\
&& \quad\qquad{}- \int_0^1 \mathbf{Q}^\top
\bigl(B' \cap A_{\dd u'}\bigr) \mathbf {I}^{-1}
\bigl(A_{u'}^c\bigr) \mathbf{Q}\bigl(B \cap
A_{u'}^c\bigr)
\\
&&\quad\qquad{} + \int_0^1 \int
_0^1 \mathbf{Q}^\top(B \cap
A_{\dd u}) \mathbf{I}^{-1}\bigl(A_u^c
\bigr) \mathbf{I}\bigl(A_{u \vee u'}^c\bigr) \mathbf
{I}^{-1}\bigl(A_{u'}^c\bigr) \mathbf{Q}
\bigl(B' \cap A_{\dd u'}\bigr).
\end{eqnarray*}
Splitting the double integral into two double integrals, one over the
region $\{u \le u'\}$ and the other over the region $\{u' < u\}$, we
see that all the integral terms cancel each other. This implies that
$W_R$ has the covariance structure of a Wiener process with time $R$.
\end{pf*}

Let $W_{n,R}$ denote the empirical version of $W_R$ in Corollary~\ref
{cor.transform},
\begin{eqnarray*}
&&
W_{n,R}\bigl([\delta,x] \times[\delta,y]\bigr)
\\
&&\qquad = \widehat{\eta}_n\bigl([\delta,x] \times[\delta,y]\bigr)
\\
&&\quad\qquad{}- \int_\delta^x \!\!\int
_\delta^y \widehat{\mathbf{q}}^\top(s,t)
\biggl( \widehat{\mathbf{I}}^{-1}(t) \int_\delta^T\!\!
\int_t^T \widehat {\mathbf{q}}
\bigl(s',t'\bigr)\,\dd\widehat{\eta}_n
\bigl(s',t'\bigr) \biggr) \,\dd R_{\widehat
{\bolds{\theta}}}(s,t).
\end{eqnarray*}
The following result will be useful for the proof of Theorem~\ref{thm.main}.

%pr7.1 #&#
\begin{prop}\label{thm.main.R}
Let $0 < \delta< \tau< T$. If assumptions \textup{B1}--\textup{B5} hold, then
\[
W_{n,R}\bigl([\delta,x] \times[\delta,y]\bigr) \Rightarrow
W_R\bigl([\delta,x] \times [\delta,y]\bigr)
\]
in $D([\delta,\tau]^2)$.
\end{prop}

\begin{pf}
Applying Skorohod's representation theorem [see, e.g., \citet{billingsley1999}, Theorem~6.7] to Theorem~\ref{etan.hat.asymp}, we
obtain a probability space that supports probabilistically equivalent
versions of $\widehat{\eta}_n$ and $\widehat{\eta}$ satisfying
\[
\| \widehat{\eta}_n - \widehat{\eta} \|_{[\delta,T]^2} \to0\qquad
\mbox{a.s.},
\]
with $\| \varphi\|_{[a,b]^2} := \sup_{(x,y) \in[a,b]^2} | \varphi
(x,y) |$. We will work on this space. Let us denote
%
%e43 #&#
%e44 #&#
%e45 #&#
%e46 #&#
\begin{eqnarray}
H(s,t) &=&  \mathbf{q}^\top(s,t)
\mathbf{I}^{-1}(t) \int_\delta^T \!\!\int
_t^T \mathbf{q}\bigl(s',t'
\bigr)\,\dd\widehat{\eta}\bigl(s',t'\bigr),
\nonumber
\\
\widehat{H}(s,t) &=& \widehat{\mathbf{q}}^\top(s,t)  \widehat{\mathbf
{I}}^{-1}(t) \int_\delta^T\!\! \int
_t^T \widehat{\mathbf{q}}\bigl(s',t'
\bigr)\,\dd\widehat{\eta}\bigl(s',t'\bigr),
\nonumber
\\[-8pt]
\label{Hs}
\\[-8pt]
\nonumber
H_n(s,t) &=& \mathbf{q}^\top(s,t)  \mathbf{I}^{-1}(t)
\int_\delta^T\!\! \int_t^T
\mathbf{q}\bigl(s',t'\bigr)\,\dd\widehat{
\eta}_n\bigl(s',t'\bigr),
\\
\nonumber
\widehat{H}_n(s,t) &=& \widehat{\mathbf{q}}^\top(s,t)
\widehat{\mathbf{I}}^{-1}(t) \int_\delta^T
\!\!\int_t^T \widehat{\mathbf {q}}
\bigl(s',t'\bigr)\,\dd\widehat{\eta}_n
\bigl(s',t'\bigr).
\end{eqnarray}
We have to show that
%
%e47 #&#
\begin{equation}
\label{H.conv}
\sup_{(x,y) \in[\delta,\tau]^2} \biggl| \int_\delta^x
\!\!\int_\delta^y \bigl( \widehat{H}_n(s,t)r_{\widehat{\bolds{\theta}}}(s,t)
- H(s,t)r_{\bolds{\theta}_0}(s,t) \bigr) \,\dd t \, \dd s \biggr| \stackrel{P}{\to} 0.
\end{equation}
For this, it suffices to prove the two statements
%
%e48 #&#
\begin{equation}
\label{ayi}
\bigl\| H ( r_{\widehat{\bolds{\theta}}} - r_{\bolds{\theta}_0} ) \bigr\|_{[\delta,\tau]^2}
\stackrel{P} {\to} 0, \qquad \bigl\| r_{\widehat
{\bolds{\theta}}} ( \widehat{H}_n - H )
\bigr\|_{[\delta,\tau]^2} \stackrel{P} {\to} 0.
\end{equation}
The first convergence in (\ref{ayi}) follows from the continuity of
$r_{\bolds{\theta}}(s,t)$ over $(s,t,\bolds{\theta}) \in
[\delta,T]^2 \times B(\bolds{\theta}_0)$ and the continuity of
$H(s,t)$ over $(s,t) \in[\delta,\tau]^2$. The second convergence in
(\ref{ayi}) follows from
%
%e49 #&#
\begin{equation}
\label{ayi3}
\| \widehat{H}_n - H \|_{[\delta,\tau]^2} \stackrel{P} {
\to} 0,
\end{equation}
since $\|r_{\widehat{\bolds{\theta}}}\|_{[\delta,\tau]^2}=O_P(1)$.
We establish (\ref{ayi3}) by proving the two statements
%
%e50 #&#
\begin{equation}
\label{esek} \| \widehat{H}_n - H_n \|_{[\delta,\tau]^2}
\stackrel{P} {\to} 0, \qquad \| H_n - H \|_{[\delta,\tau]^2}
\stackrel{P} {\to} 0.
\end{equation}
Consider the second statement in (\ref{esek}). Its left-hand side is
equal to
\[
\sup_{(s,t) \in[\delta,\tau]^2}\biggl| \mathbf{q}^\top(s,t)  \mathbf
{I}^{-1}(t) \int_\delta^T \!\!\int
_t^T \mathbf{q}\bigl(s',t'
\bigr)\,\dd\Delta _n\bigl(s',t'\bigr) \biggr|,
\]
with $\Delta_n = \widehat{\eta}_n - \widehat{\eta}$. The vector
function $\mathbf{q}^\top(s,t) \mathbf{I}^{-1}(t)$ is bounded on $(s,t)
\in[\delta,\tau]^2$, by continuity. So it will suffice to show
%
%e51 #&#
\begin{equation}
\label{esek3}
\sup_{t \in[\delta,\tau]} \biggl| \int_\delta^T
\!\!\int_t^T q_i\bigl(s',t'
\bigr)\, \dd\Delta_n\bigl(s',t'\bigr) \biggr|
\stackrel{P} {\to} 0, \qquad  i=1,\ldots,6+d.
\end{equation}
The double integral inside the absolute value bars can be rewritten,
using integration by parts [see \citet{hildebrandt1963}, Section~III.8], as follows:
%
%e52 #&#
%e53 #&#
%e54 #&#
%e55 #&#
\begin{eqnarray*}
&& q_i(T,T)\Delta_n(T,T) - q_i(T,t)
\Delta_n(T,t) - q_i(\delta ,T)\Delta_n(
\delta,T) + q_i(\delta,t)\Delta_n(\delta,t)
\\
&&\qquad{}- \int_\delta^T \Delta_n
\bigl(s',T\bigr) \, \dd q_i\bigl(s',T
\bigr) + \int_\delta^T \Delta_n
\bigl(s',t\bigr) \, \dd q_i\bigl(s',t
\bigr)
\\
&&\qquad{}- \int_t^T \Delta_n
\bigl(T,t'\bigr) \, \dd q_i\bigl(T,t'
\bigr) + \int_t^T \Delta_n\bigl(
\delta ,t'\bigr) \, \dd q_i\bigl(\delta,t'
\bigr)
\\
&&\qquad{}+ \int_\delta^T\!\! \int_t^T
\Delta_n\bigl(s',t'\bigr) \, \dd
q_i\bigl(s',t'\bigr).
\end{eqnarray*}
Each of the first four terms is bounded in absolute value by $\|q_i\|
_{[\delta,T]^2} \cdot\|\Delta_n\|_{[\delta,T]^2}$, where the first
factor is finite by continuity and the second factor vanishes in
probability. Moreover, each integral term is bounded in absolute value
by $\| q_i \|_{\mathrm{HK}} \|\Delta_n\|_{[\delta,T]^2}$, which also
vanishes in probability because $\|q_i\|_{\mathrm{HK}} < \infty$, by
virtue of the assumptions $\|\rho_j(x,y,\bolds{\theta}_0)\|_{\mathrm
{HK}} < \infty$ for $j=1,\ldots,2+d$, and Proposition~1 of \citet{blumlingertichy1989}. Hence (\ref{esek3}) follows, and the second
convergence in (\ref{esek}) is established.

It remains to prove the first convergence in (\ref{esek}). By virtue of
the second convergence there, and an analogous\ result for $\widehat
{H}_n$ and $\widehat{H}$, it will suffice to prove $\| \widehat{H} - H
\|_{[\delta,\tau]^2} \stackrel{P}{\to} 0$. Note that
%
%e56 #&#
%e57 #&#
%e58 #&#
\begin{eqnarray}
&& \bigl| \widehat{H}(s,t) - H(s,t) \bigr|\nonumber
\\
\label{H.decomp}
&& \qquad \le \bigl| \widehat{\mathbf{q}}^\top(s,t) \widehat{\mathbf{I}}^{-1}(t)
- \mathbf{q}^\top(s,t) \mathbf{I}^{-1}(t) \bigr| \cdot \biggl| \int
_\delta^T\!\! \int_t^T
\mathbf{q}\bigl(s',t'\bigr) \, \dd \widehat{\eta}
\bigl(s',t'\bigr) \biggr|
\\
&& \quad\qquad{}+ \bigl| \widehat{\mathbf{q}}^\top(s,t) \widehat{
\mathbf{I}}^{-1}(t) \bigr| \cdot \biggl| \int_\delta^T
\!\!\int_t^T \bigl( \widehat{\mathbf {q}}
\bigl(s',t'\bigr) - \mathbf{q}\bigl(s',t'
\bigr) \bigr)\, \dd\widehat{\eta}\bigl(s',t'\bigr) \biggr|,\nonumber
\end{eqnarray}
where $| \, \cdot\, |$ should be interpreted component-wise.

Let us write $\mathbf{q}(s,t,z_1,z_2,w_1,\ldots,w_d)$ to denote the
vector $\mathbf{q}(s,t)$ with the values of $\gamma_1$ and $\gamma_2$
replaced by variables $z_1$ and $z_2$, and the values $\theta_{01},
\ldots, \theta_{0d}$ replaced by variables $w_1,\ldots,w_d$. Then
$\mathbf{q}(s,t) = \mathbf{q}(s,t,\gamma_1, \gamma_2, \theta_{01},
\ldots, \theta_{0d})$ and $\widehat{\mathbf{q}}(s,t) = \mathbf
{q}(s,t,\widehat{\gamma}_1, \widehat{\gamma}_2, \widehat{\theta}_1,
\ldots, \widehat{\theta}_d)$.

Now consider the first term on the right-hand side of (\ref{H.decomp}).
Since the vector $\mathbf{q}(s,t,z_1,z_2,w_1,\ldots,w_d)$ is continuous\vspace*{1.5pt}
over $[\delta,\tau]^2 \times\mathbb{R}^2 \times B(\bolds{\theta
}_0)$, we have that $|\widehat{\mathbf{q}}^\top(s,t) \widehat{\mathbf
{I}}^{-1}(t) - \mathbf{q}^\top(s,t) \mathbf{I}^{-1}(t)|$ is $ o_P(1)$
uniformly over $(s,t) \in[\delta,\tau]^2$. Moreover, an integration by
parts argument as above yields that
\[
\biggl| \int_\delta^T\!\! \int_t^T
q_i\bigl(s',t'\bigr) \, \dd\widehat{
\eta}\bigl(s',t'\bigr) \biggr| \le\| \widehat{\eta}
\|_{[\delta,T]^2} \cdot\bigl(4\| q_i \| _{[\delta,T]^2} + 5\|
q_i \|_{\mathrm{HK}}\bigr)
\]
for $1 \le i \le6+d$, where the right-hand side is $O_P(1)$. We
conclude that the first term on the right-hand side of (\ref{H.decomp})
is $o_P(1)$ uniformly over $(s,t) \in[\delta,\tau]^2$.

Next,\vspace*{1.5pt} consider the second term on the right-hand side of (\ref
{H.decomp}). It follows from the discussion above that the vector $
| \widehat{\mathbf{q}}^\top(s,t)\widehat{\mathbf{I}}^{-1}(t)  |$ is
$O_P(1)$ uniformly over $(s,t) \in[\delta,\tau]^2$, so it will suffice
to show that
%
%e59 #&#
\begin{equation}
\label{kurt}
\sup_{t \in[\delta,\tau]} \biggl| \int_\delta^T
\!\!\int_t^T \Delta q_i
\bigl(s',t'\bigr) \,\dd\widehat{\eta}
\bigl(s',t'\bigr) \biggr| \stackrel{P} {\to} 0,
\end{equation}
with $\Delta q_i = \widehat{q}_i - q_i$, for $i=1,\ldots,6+d$. Once
again, an integration by parts argument shows that the left-hand side
of (\ref{kurt}) is bounded from above by
\[
\| \widehat{\eta} \|_{[\delta,T]^2} \cdot\bigl(4\| \Delta q_i
\|_{[\delta
,T]^2} + 5\| \Delta q_i \|_{\mathrm{HK}}\bigr),
\]
where $\| \widehat{\eta} \|_{[\delta,T]^2} < \infty$ a.s. and $\|
\Delta q_i \|_{[\delta,T]^2} = o_P(1)$ by continuity.\vspace*{1pt} It remains to
establish $\|\Delta q_i\|_{\mathrm{HK}} =o_P(1)$. For $i=7,\ldots,6+d$,
this follows directly from assumption B5. For $i=1$, we have
\begin{eqnarray*}
\|\Delta q_1 \|_{\mathrm{HK}} &=& \bigl\| f_1(x,\widehat{
\gamma}_1) \rho_1(x,y,\widehat{\bolds{\theta }}) -
f_1(x,\gamma_1) \rho_1(x,y,\bolds{
\theta}_0) + \Delta f'_1(x)
\bigr\|_{\mathrm{HK}}
\\
&\le& \bigl\| \Delta f_1(x) \rho_1(x,y,\bolds{
\theta}_0) \bigr\|_{\mathrm{HK}} + \bigl\| f_1(x,\widehat{
\gamma}_1) \Delta\rho_1(x,y) \bigr\|_{\mathrm{HK}} +
2V^{(1)}\bigl(\Delta f'_1\bigr).
\end{eqnarray*}
Using Proposition~1 of \citet{blumlingertichy1989}, differentiability
properties of $f_1$, $f'_1$ on $[\delta,T]$ and assumption B5, each
term on the right-hand side can be shown to be $o_P(1)$. The cases
$i=2,\ldots,6$ are similar. Thus (\ref{kurt}) follows.
\end{pf}

\begin{pf*}{Proof of Theorem~\protect\ref{thm.main}}
Note that we have
\begin{eqnarray*}
W_n\bigl([\delta,x] \times[\delta,y]\bigr) &=&  \int
_\delta^x \!\!\int_\delta^y
\sigma (s,t,\widehat{\bolds{\theta}}) \, \dd W_{n,R}\bigl([\delta,s]
\times [\delta,t]\bigr),
\\
W\bigl([\delta,x] \times[\delta,y]\bigr) &=& \int_\delta^x
\!\!\int_\delta^y \sigma (s,t,\bolds{
\theta}_0) \, \dd W_{R}\bigl([\delta,s] \times[\delta,t]
\bigr).
\end{eqnarray*}
Now, by Proposition~\ref{thm.main.R} and Skorohod's representation
theorem, there exists a probability space supporting versions of
$W_{n,R}$ and $W_R$ which satisfy
\begin{eqnarray*}
&& \sup_{(x,y) \in[\delta,\tau]^2} \bigl| W_{n,R}\bigl([\delta,x]
\times[\delta,y]\bigr) - W_R\bigl([\delta,x] \times[\delta,y]\bigr)
\bigr|
\\
&&\qquad =: \sup_{(x,y) \in[\delta,\tau]^2}\bigl| D_n(x,y) \bigr| \to0 \qquad  \mbox{a.s.}
\end{eqnarray*}
We work with this probability space. We have
%
%e60 #&#
%e61 #&#
\begin{eqnarray}
&& \bigl| W_n\bigl([\delta,x] \times[\delta,y]\bigr)
- W\bigl([\delta,x] \times [\delta,y]\bigr) \bigr|
\nonumber
\\[-8pt]
\label{hayvan}
\\[-8pt]
\nonumber
&&\qquad \le \biggl| \int_\delta^x \!\!\int_\delta^y
\Delta\sigma(s,t) \, \dd W_R\bigl([\delta,s] \times[\delta,t]\bigr)
\biggr| + \biggl| \int_\delta^x \!\!\int_\delta^y
\sigma(s,t,\widehat{\bolds {\theta}}) \, \dd D_n(s,t) \biggr|.\hspace*{-20pt}
\end{eqnarray}
Applying integration by parts as in the proof of Proposition~\ref
{thm.main.R}, we see that the first term on the right-hand side of (\ref
{hayvan}) is bounded by
%
%e62 #&#
\begin{equation}
\label{yunus} \sup_{(s,t) \in[\delta,\tau]^2} \bigl|W_R\bigl([\delta,s]
\times[\delta,t]\bigr)\bigr| \cdot\bigl(4 \| \Delta\sigma\|_{[\delta,\tau]^2} + 5 \| \Delta
\sigma\| _{\mathrm{HK}}\bigr).
\end{equation}
Since $W_R$ is a.s. bounded on $[\delta,\tau]^2$, $\| \Delta\sigma\|
_{[\delta,\tau]^2} = o_P(1)$ by continuity, and $\|\Delta\sigma\|_{\mathrm{HK}} = o_P(1)$ by assumption B5, (\ref{yunus}) vanishes in probability.
Similarly, the second term on the right-hand side of (\ref{hayvan}) is
bounded by
\[
\|D_n\|_{[\delta,\tau]^2} \cdot\bigl(4 \bigl\| \sigma(\cdot,\cdot,\widehat
{\bolds{\theta}}) \bigr\|_{[\delta,\tau]^2} + 5 \bigl\| \sigma(\cdot,\cdot ,\widehat{\bolds{
\theta}}) \bigr\|_{\mathrm{HK}}\bigr),
\]
which also vanishes in probability since $\|D_n\|_{[\delta,\tau]^2} =
o_P(1)$ and the two summands in the parentheses are $O_P(1)$. Thus the
left-hand side of (\ref{hayvan}) is $o_P(1)$ uniformly over $(s,t) \in
[\delta,\tau]^2$.
\end{pf*}

\section*{Acknowledgements}
We are very grateful to an Associate Editor and three referees for many
insightful comments and suggestions that led to this improved version
of the manuscript. We are also grateful to the participants of the 2013
Extreme Value Analysis conference in Shanghai and the Fourth Wellington
Workshop in Probability and Mathematical Statistics (WWPMS4) for their feedback.

\begin{supplement}[id=suppA]
\stitle{Supplement to ``Asymptotically distribution-free
goodness-of-fit testing for tail copulas''}
\slink[doi]{10.1214/14-AOS1304SUPP} %[doi,text={...}] - jei reikia
%suskaldyti doi
\sdatatype{.pdf}
\sfilename{aos1304\_supp.pdf}
\sdescription{We provide a proof of Theorem~\ref{etan.asymp} as well
as details about the Monte Carlo simulations of Section~\ref{sec.sim}.}
\end{supplement}

%
% imsref loaded by daiva.urboniene, 2015-02-23 10:32:34


\printaddresses
\end{document}